\documentclass[a4paper]{amsart}
\usepackage{amsmath}
\usepackage{amssymb}
\usepackage[all]{xy}
\xyoption{matrix}
\xyoption{arrow}
\begin{document}

\newtheorem{theo}{Theorem}[section]
\newcommand{\bthm}{\begin{theo}}
\newcommand{\ethm}{\end{theo}}
\newtheorem{coro}[theo]{Corollary}
\newcommand{\bcor}{\begin{coro}}
\newcommand{\ecor}{\end{coro}}
\newtheorem{prop}[theo]{Proposition}
\newcommand{\bprop}{\begin{prop}}
\newcommand{\eprop}{\end{prop}}
\newtheorem{claim}[theo]{Lemma}
\newcommand{\blemma}{\begin{claim}}
\newcommand{\elemma}{\end{claim}}
\newtheorem{conjec}[theo]{Conjecture}
\newcommand{\bconj}{\begin{conjec}}
\newcommand{\econj}{\end{conjec}}

\theoremstyle{definition}
\newtheorem{defo}[theo]{Definition}
\newcommand{\bdfn}{\begin{defo}}
\newcommand{\edfn}{\end{defo}}
\newtheorem{rmk}[theo]{Remark}
\newcommand{\brk}{\begin{rmk}}
\newcommand{\erk}{\end{rmk}}
\newtheorem{exple}[theo]{Example}
\newcommand{\bex}{\begin{exple}}
\newcommand{\eex}{\end{exple}}

\newcommand{\note}{\mbox{\textbf{Note: \hspace{0.2cm}}}}


\title[Companion Bases for Cluster-tilted Algebras]{Companion Bases for Cluster-tilted Algebras}

\author[Parsons]{Mark James Parsons}
\address{Berlin \\
Germany
}
\email{markjamesparsons@googlemail.com}


\begin{abstract}
Motivated by work of Barot, Geiss and Zelevinsky, we study a collection of $\mathbb{Z}$-bases (which we call \textit{companion bases}) of the integral root lattice of a root system of simply-laced Dynkin type. Each companion basis is associated with the quiver of a cluster-tilted algebra of the corresponding type.

In type $A$, we establish that the dimension vectors of the finitely generated indecomposable modules over a cluster-tilted algebra may be obtained, up to sign, by expanding the positive roots in terms of any companion basis for the quiver of that algebra. This generalises part of Gabriel's Theorem.

Also, we describe the relationship between different companion bases for the same quiver and show how to mutate a companion basis for a quiver to produce a companion basis for a mutation of that quiver.
\end{abstract}

\date{11 October 2011}


\thanks{This work was supported by the
Engineering and Physical Sciences Research Council}

\maketitle


\section{Introduction}
\label{introduction}

Cluster algebras were introduced by Fomin and Zelevinsky \cite{FZ1} in order to better understand the dual canonical basis of the quantised enveloping algebra of a finite dimensional semisimple Lie algebra. In \cite{MRZ}, a link between cluster algebras and representations of quivers was established. This subsequently led to the introduction of cluster categories in \cite{BMRRT}, which were intended to give a categorical model of cluster algebras. (Note that independently of \cite{BMRRT}, a geometric definition of cluster categories of Dynkin type $A$ was given in \cite{CCS1}.) A key development in the study of cluster categories was the creation of a generalised version of APR-tilting theory (see \cite{APR}), known as cluster-tilting theory. In this theory, a key role is played by the cluster-tilted algebras, as introduced in \cite{BMR1}.

In \cite{BMR2}, an important link between cluster-tilted algebras and cluster algebras was established. Each cluster algebra is associated with an equivalence class of skew-symmetrizable integer matrices, known as the exchange matrices of that cluster algebra. In particular, the exchange matrices of a cluster algebra of simply-laced Dynkin type are skew-symmetric, and can therefore be represented as quivers. It was shown independently in \cite{BMR2} and \cite{CCS2} that the quivers of the cluster-tilted algebras of a given simply-laced Dynkin type are precisely the quivers of the exchange matrices of the cluster algebra of that type. We present the background material we require on cluster algebras, cluster categories and cluster-tilted algebras in Section~\ref{Cluster Algebras and Cluster Categories}.

In addressing the problem of recognising the cluster algebras of finite type, the paper \cite{BGZ} considered a class of matrices closely related to the Cartan matrices, known as the positive quasi-Cartan matrices. We recall the main results of \cite{BGZ}, including their results on positive quasi-Cartan matrices, in Section~\ref{Positive Quasi-Cartan Matrices}. Given a cluster-tilted algebra of simply-laced Dynkin type, it follows from the main result of \cite{BGZ} that the exchange matrix associated to its quiver must have a positive quasi-Cartan companion. (That is, there exists some positive quasi-Cartan matrix for which the absolute values of the off-diagonal entries match the absolute values of the off-diagonal entries of the exchange matrix.) Of key motivational importance is the classification result of the positive quasi-Cartan matrices, which tells us that this positive quasi-Cartan companion arises as the matrix of inner products associated to some $\mathbb{Z}$-basis of roots of the integral root lattice of the corresponding root system of simply-laced Dynkin type. We call such a $\mathbb{Z}$-basis of roots of the integral root lattice a companion basis giving rise to that particular positive quasi-Cartan companion.

Companion bases are formally defined in Definition~\ref{companion basis}. The remainder of Section~\ref{companion bases} examines the relationship between different companion bases which give rise to the same positive quasi-Cartan companion of an exchange matrix. A complete description of this relationship is obtained in Theorem~\ref{description of beta sets giving rise to A}. Furthermore, as a simple consequence of this we are able to describe, in Corollary~\ref{description of beta sets}, the relationship between any two companion bases giving rise to different positive quasi-Cartan companions of the same exchange matrix.

Given a cluster-tilted algebra of simply-laced Dynkin type, to any companion basis giving rise to a positive quasi-Cartan companion of the exchange matrix associated to its quiver, we associate a collection of vectors. These vectors are obtained by expressing each positive root of the corresponding root system in terms of the elements of the companion basis, and taking the absolute values of the coefficients. Our main result, Theorem~\ref{main result}, shows that in the Dynkin type $A$ case, these vectors are precisely the dimension vectors of the finitely generated indecomposable modules over the given cluster-tilted algebra. Section~\ref{A Generalisation of Gabriel's Theorem} is devoted to proving this result which can be regarded as a generalisation, in the Dynkin type $A$ case, of part of Gabriel's Theorem \cite{Gab}. The proof takes advantange of a well-known description of the quivers of the cluster-tilted algebras of Dynkin type $A$ in terms of triangulations of regular polygons (see \cite{CCS1,CCS2}).

We conjecture that our main result may be extended to all of the simply-laced Dynkin cases (Conjecture~\ref{main conjecture}). This provides our motivation for the material presented in Sections \ref{Companion Basis Mutation and Mutation Maps} and \ref{Towards a Description of the Mutation Map}.

Each companion basis has an associated exchange matrix, and hence may be associated to the quiver of a cluster-tilted algebra of simply-laced Dynkin type. The main result of Section~\ref{Companion Basis Mutation and Mutation Maps} is Theorem~\ref{beta set mutation theorem} which establishes a companion basis mutation procedure that, given a companion basis for the quiver of a cluster-tilted algebra of simply-laced Dynkin type, produces a companion basis for any mutation of that quiver.

We show how companion basis mutation induces a map from the collection of vectors associated to a companion basis for the quiver of a cluster-tilted algebra of simply-laced Dynkin type to the collection of vectors associated to the mutated companion basis for a mutation of the initial quiver. Moreover, Theorem~\ref{mutation maps all the same} establishes that this induced map is independent of the choice of companion basis for the initial quiver. Finally, Section~\ref{Towards a Description of the Mutation Map} gives a description of the induced mutation maps in the Dynkin type $A$ case. As an interesting consequence, this description shows how the dimension vectors of the finitely generated indecomposable modules over a cluster-tilted algebra of Dynkin type $A$ may be used to write down the dimension vectors of the finitely generated indecomposable modules over a cluster-tilted algebra obtained from the former algebra by cluster-tilting once.

We remark that independent work of Ringel \cite{Rin2} gives an alternative approach for computing the dimension vectors of the finitely generated indecomposable modules over a cluster-tilted algebra, whenever the associated cluster-tilting object corresponds to a preprojective tilting module.


\section{Cluster Algebras and Cluster Categories}
\label{Cluster Algebras and Cluster Categories}

Cluster algebras were introduced by Fomin and Zelevinsky \cite{FZ1} in 2001. A cluster algebra is a subring of the field of rational functions in $n$ indeterminates, generated by cluster variables. These cluster variables are obtained via a process of mutation, starting from some `initial seed'. Cluster algebras having only finitely many cluster variables are known as cluster algebras of finite type. One of the early key results in the development of the theory of cluster algebras was the classification of the cluster algebras of finite type, given in \cite{FZ2}.

We start by recalling the definition of a cluster algebra and the classification of the cluster algebras of finite type, but first, we need a preliminary definition.

\bdfn
\label{skew symmetrizable}
A square integer matrix $B$ is said to be \textit{skew-symmetrizable} (resp.\ \textit{symmetrizable}) if there is some diagonal matrix $D$ with positive integer diagonal entries such that $DB$ is skew-symmetric (resp.\ symmetric).
\edfn

We can now give the definition of a cluster algebra (without coefficients).

Let $n \in \mathbb{N}$ and let $\mathbb{F} = \mathbb{Q} (u_1, \ldots, u_n)$ be the field of rational functions in $n$ indeterminates. Let $\mathbf{x} = \{ x_1, \ldots, x_n \} \subseteq \mathbb{F}$ be a free generating set for $\mathbb{F}$ over $\mathbb{Q}$, and let $B = (b_{xy})_{x,y \in \mathbf{x}}$ be an $n \times n$ skew-symmetrizable integer matrix with rows and columns indexed by the entries of $\mathbf{x}$. We call the pair $(\mathbf{x},B)$ a seed in $\mathbb{F}$, and we obtain more seeds from this `initial seed' via a mutation process. Let $z \in \mathbf{x}$. We obtain a new free generating set $\mathbf{x}' = (\mathbf{x} \setminus \{ z \} ) \cup \{ z' \}$ for $\mathbb{F}$ over $\mathbb{Q}$, where $z'$ is obtained using the exchange relation
$$
zz' = \prod_{x \in \mathbf{x}, \, b_{xz}>0} x^{b_{xz}} + \prod_{x \in \mathbf{x}, \, b_{xz}<0} x^{-b_{xz}}.
$$
Similarly, we obtain a skew-symmetrizable matrix (see \cite[Proposition 4.5]{FZ1}) $B'$ from $B$ with entries given by
\begin{eqnarray*}
\label{exchange matrix}
b'_{xy} = \left \{ \begin{array}{ll}
-b_{xy} & \textrm{if } x=z \textrm{ or } y=z,\\
b_{xy}+ \frac{|b_{xz}|b_{zy}+b_{xz}|b_{zy}|}{2} & \textrm{otherwise}.
\end{array} \right.
\end{eqnarray*}
The row and column labelled $z$ in $B$ are relabelled $z'$ in $B'$. The pair $(\mathbf{x}' , B')$ then form a seed which we call the mutation of $(\mathbf{x} , B)$ in the direction $z$. It can easily be checked that by mutation of $(\mathbf{x}' , B')$ in the direction $z'$, we recover the seed $(\mathbf{x} , B)$. By iterated mutations of the initial seed $(\mathbf{x} , B)$ in all directions, we obtain a set of seeds $\mathcal{S}$. The free generating sets appearing in these seeds are called clusters, and their elements are called cluster variables. The matrices appearing in these seeds are called exchange matrices. The set of all cluster variables is the union of all of the free generating sets appearing in seeds in $\mathcal{S}$, and is denoted by $\chi$. The \textit{cluster algebra} $\mathcal{A} = \mathcal{A} (\mathbf{x} , B)$ is then defined to be the $\mathbb{Z}$-subalgebra (subring) of $\mathbb{F}$ generated by $\chi$.

Two cluster algebras $\mathcal{A}'$ and $\mathcal{A}''$, contained in fields of rational functions $\mathbb{F}'$ and $\mathbb{F}''$ respectively, are said to be \textit{isomorphic} (as cluster algebras) if there is a $\mathbb{Z}$-algebra isomorphism $\iota \colon \mathbb{F}' \rightarrow \mathbb{F}''$ taking some seed $(\mathbf{y} , C)$ of $\mathcal{A}'$ to a seed $(\iota (\mathbf{y}) , C)$ of $\mathcal{A}''$. (We note that the terminology ``strongly isomorphic'' is used in \cite{FZ2}.)

Up to isomorphism of cluster algebras, $\mathcal{A} (\mathbf{x} , B)$ does not depend on the choice of free generating set $\mathbf{x}$ for $\mathbb{F}$, and so we can denote this cluster algebra by $\mathcal{A} (B)$. In fact, we can go further than this. It is clear that the mutation of skew-symmetrizable matrices outlined above gives rise to an equivalence relation on the set of $n \times n$ skew-symmetrizable integer matrices. Up to isomorphism of cluster algebras, the cluster algebra $\mathcal{A} (B)$ only depends on the mutation equivalence class of $B$.

We have the following definition from \cite{FZ2}.

\bdfn
\label{cluster algebra finite type}
A cluster algebra is said to be of \textit{finite type} if it has only finitely many cluster variables.
\edfn

Before stating the result classifying the cluster algebras of finite type, we need to introduce some terminology.

\bdfn
\label{cartan counterpart}
Let $B = (b_{ij})$ be an $n \times n$ integer matrix. Then the \textit{Cartan counterpart} of $B$ is defined to be the matrix $A(B) = (a_{ij})$ given by $a_{ii} = 2$ for all $1 \leq i \leq n$, and $a_{ij} = -|b_{ij}|$ otherwise.
\edfn

\note
It is clear that the Cartan counterpart of a skew-symmetrizable matrix must be a symmetrizable matrix.

We can now give the classification result of the cluster algebras of finite type. The result was originally given in \cite{FZ2}, but we give the version as stated in \cite{BGZ}.

\bthm
\label{finite type classification}
Let $\mathcal{F}$ be a mutation equivalence class of skew-symmetrizable matrices. Then the following are equivalent:\\
(i) The cluster algebra associated to $\mathcal{F}$ is of finite type.\\
(ii) There is some matrix $B \in \mathcal{F}$ such that the Cartan counterpart of $B$ is a Cartan matrix of finite type.\\
(iii) For every matrix $B = (b_{ij}) \in \mathcal{F}$, $|b_{ij} b_{ji}| \leq 3$ for all $i \neq j$.

Furthermore, the Cartan-Killing type of the Cartan matrix in (ii) is uniquely determined by the equivalence class $\mathcal{F}$ (and is referred to as being the type of the cluster algebra associated to $\mathcal{F}$).
\ethm

Since the Cartan matrices of finite type can be represented by Dynkin diagrams, we can consider this result as giving as giving a classification of the cluster algebras of finite type by Dynkin diagrams.

If $B = (b_{ij})$ is an $n \times n$ skew-symmetric integer matrix, then we can associate to $B$ a quiver $\Gamma (B)$ with vertices corresponding to the rows and columns of $B$, and $b_{ij}$ arrows from the vertex $i$ to the vertex $j$ whenever $b_{ij}>0$. It is easy to show that mutation in any direction of a skew-symmetric matrix results in a skew-symmetric matrix, and hence the mutation equivalence class of $B$ must consist entirely of skew-symmetric matrices. Therefore, we can associate a quiver to each matrix belonging to the mutation equivalence class of $B$.

In view of Theorem~\ref{finite type classification}, we see that the cluster algebras of Dynkin types $A$, $D$ and $E$ are those associated with equivalence classes of skew-symmetric matrices which contain a matrix whose associated quiver is a Dynkin quiver of the same type.

\brk
\label{matrix entries remark}
Suppose $C = (c_{ij})$ is a matrix appearing in a seed of a cluster algebra of (simply-laced) Dynkin type $A$, $D$ or $E$. Then, we have that $C$ is a skew-symmetric integer matrix, and also, from Theorem~\ref{finite type classification}, we see that $C$ must satisfy $|c_{ij} c_{ji}| \leq 3$ for all $i \neq j$. It follows that the entries of $C$ must all belong to the set $\{ 0, \pm 1 \}$.
\erk

We now turn our attention to cluster categories. Following \cite{BMRRT}, we state the definition of a cluster category of simply-laced Dynkin type and briefly consider some aspects of the relationship between such a cluster category and the corresponding cluster algebra.

Let $k$ be an algebraically closed field and let $Q$ be a simply-laced quiver of Dynkin type with underlying graph $\Delta$. Let $\mathcal{D} = \mathcal{D}^b (kQ \mathrm{-mod})$ be the bounded derived category of the category of finitely generated left $kQ$-modules with shift functor $[1]$. Also, let $\tau$ be the AR-translation in $\mathcal{D}$, and define $F = \tau^{-1} [1]$. The \textit{cluster category} $\mathcal{C} = \mathcal{C} (kQ)$ is then defined to be the factor category
$$
\mathcal{C} = \frac{\mathcal{D}^b (kQ \mathrm{-mod})}{F}
$$
whose objects are the objects of $\mathcal{D}$, and where morphisms are given by
$$
\mathrm{Hom}_{\mathcal{C}} (X,Y) = \bigoplus_{i \in \mathbb{Z}} \mathrm{Hom}_{\mathcal{D}} (X,F^iY)
$$
for objects $X,Y$ in $\mathcal{C}$.

Let $\Phi$ be the root system of Dynkin type $\Delta$. Suppose that $\Pi = \{ \alpha_1 , \ldots , \alpha_n \} \subseteq \Phi$ is a simple system of $\Phi$, and that $\Phi^+$, $\Phi_{\geq -1} = \Phi^+ \cup (-\Pi)$ are respectively the corresponding sets of positive and almost positive roots. In \cite{BMRRT}, it is observed that there is a one-to-one correspondence between the set of isomorphism classes of indecomposable objects in $\mathcal{C}$ and the set of almost positive roots $\Phi_{\geq -1}$. Indeed, by identifying an indecomposable left $kQ$-module $M$ with the stalk complex where $M$ appears in degree zero, we may consider $M$ as an object in $\mathcal{C}$. These objects together with the objects of the form $P[1]$ for $P$ a finitely generated indecomposable left $kQ$-module are, up to isomorphism, precisely the indecomposable objects in $\mathcal{C}$. Gabriel's Theorem \cite{Gab} provides a correspondence between the indecomposable $kQ$-modules and the set of positive roots $\Phi^+$. The picture is completed by associating, for each vertex $i$ of $Q$, the negative simple root $-\alpha_i$ to the object of the form $P_i [1]$ where $P_i$ is the indecomposable projective $kQ$-module corresponding to $i$.

Let $\mathcal{A}$ be the cluster algebra of Dynkin type $\Delta$. The result \cite[Theorem 1.9]{FZ2} exhibited a unique bijection between the cluster variables of $\mathcal{A}$ and the almost positive roots in $\Phi$, obtained by expressing each cluster variable in terms of a specially chosen initial seed. The denominator of each such expression is a monomial in the cluster variables of the initial seed, and the corresponding almost positive root is obtained as a linear combination of simple roots with coefficients given by the exponents appearing in this monomial.

This demonstrates a link between the cluster category $\mathcal{C}$ and the cluster algebra $\mathcal{A}$, as combining the above provides a one-to-one correspondence between the set of isomorphism classes of indecomposable objects in $\mathcal{C}$ and the set of cluster variables in $\mathcal{A}$.

\bdfn
\label{cluster-tilting object}
An object $T$ in $\mathcal{C}$ is said to be \textit{cluster-tilting} provided for any object $X$ in $\mathcal{C}$, we have $\mathrm{Ext}^1 _{\mathcal{C}} (T,X) = 0$ if and only if $X$ lies in the additive subcategory $\mathrm{add} (T)$ of $\mathcal{C}$ generated by $T$.

A cluster-tilting object is said to be \textit{basic} if all of its direct summands are non-isomorphic.
\edfn

In our context, the cluster-tilting objects in $\mathcal{C}$ are precisely the maximal rigid objects in $\mathcal{C}$. That is, the objects $T$ in $\mathcal{C}$ for which $\mathrm{Ext}^1 _{\mathcal{C}} (T,T) = 0$ and $T$ is maximal with this property. (See \cite[Section 2.3]{Rei}, for example.) Also, we note from \cite[Theorem 3.3]{BMRRT} that the number of indecomposable direct summands of a basic cluster-tilting object in $\mathcal{C}$ is equal to the number of simple $kQ$-modules (i.e.\ the number of vertices of $Q$).

Associated to the cluster-tilting objects in $\mathcal{C}$ are the cluster-tilted algebras.

\bdfn
\label{cluster-tilted algebra}
Let $T$ be a cluster-tilting object in $\mathcal{C}$. The cluster-tilted algebra (of Dynkin type $\Delta$) associated to $T$ is the algebra $\mathrm{End}_{\mathcal{C}} (T)^{\mathrm{op}}$.
\edfn

The link between $\mathcal{C}$ and $\mathcal{A}$ noted above, can be seen to be even deeper by making use of the compatibility degree $( \; || \; ) \colon \Phi_{\geq -1} \times \Phi_{\geq -1} \rightarrow \mathbb{N} \cup \{ 0  \}$, introduced in \cite{FZ3}. We omit the definition of the compatibility degree here as it suffices for our purposes to note two of its properties. For this, we first need to mention that two almost positive roots are said to be compatible if their compatibility degree is 0. The first property is that under the bijection between the cluster variables of $\mathcal{A}$ and almost positive roots $\Phi_{\geq -1}$, maximal pairwise compatible subsets of $\Phi_{\geq -1}$ correspond to clusters in $\mathcal{A}$. We now assume $Q$ to be an alternating quiver (but see \cite[Section 3]{Zhu}). In this case, if $\alpha , \beta \in \Phi_{\geq -1}$ and $M_{\alpha} , M_{\beta}$ are the respective corresponding (isomorphism classes of) indecomposable objects in $\mathcal{C}$, then we have from \cite[Corollary 4.3]{BMRRT} that $( \alpha || \beta ) = \mathrm{dim} \, \mathrm{Ext}^1 _\mathcal{C} (M_{\alpha} , M_{\beta})$. This establishes that the bijection between indecomposable objects in $\mathcal{C}$ and cluster variables in $\mathcal{A}$ induces a bijection between basic cluster-tilting objects and clusters, as it shows that the set of indecomposable direct summands of a basic cluster-tilting object in $\mathcal{C}$ must correspond to a maximal compatible set of almost positive roots. Furthermore, since each seed in $\mathcal{A}$ is uniquely determined by its cluster \cite[Theorem 1.12]{FZ2}, we see that there is a one-to-one correspondence between the set of basic cluster-tilting objects in $\mathcal{C}$ and the set of seeds of $\mathcal{A}$.

Finally, we state a result of \cite[Section 6]{BMR2}, \cite[Theorem 3.1]{CCS2} which says that the quiver of the cluster-tilted algebra associated to a basic cluster-tilting object is precisely the quiver associated to the exchange matrix appearing in the corresponding seed.

\bthm
\label{quiver of cta is quiver of seed}
Let $T$ be a basic cluster-tilting object in $\mathcal{C}$ and suppose that the seed of the cluster algebra $\mathcal{A}$ corresponding to $T$ is $(\mathbf{x},B)$. Let $\Lambda$ be the cluster-tilted algebra $\Lambda = \mathrm{End}_{\mathcal{C}} (T)^{\mathrm{op}}$, and suppose that $Q_{\Lambda}$ is the quiver of $\Lambda$. Then, $Q_{\Lambda}$ is identical to the quiver $\Gamma (B)$ associated to the skew-symmetric matrix $B$ (identifying indecomposable objects in $\mathcal{C}$ with the corresponding cluster variables in $\mathcal{A}$).
\ethm


\section{Positive Quasi-Cartan Matrices}
\label{Positive Quasi-Cartan Matrices}

Given a skew-symmetrizable matrix $B$, the classification theorem for the cluster algebras of finite type provides two conditions for checking whether or not the cluster algebra $\mathcal{A} = \mathcal{A} (B)$ is of finite type. The paper \cite{BGZ} highlights that both of these conditions can be difficult to check in general. The focus of that paper is to solve this problem by giving a method for determining whether or not the cluster algebra $\mathcal{A} (B)$ is of finite type, based solely on consideration of the matrix $B$ itself. The solution presented in \cite{BGZ} makes use of positive quasi-Cartan matrices.

In the first part of this section, we state the main result of \cite{BGZ}, and also two further results from \cite{BGZ} which highlight its usefulness. We then focus on the positive quasi-Cartan matrices. In particular, we state a result classifying the equivalence classes of these matrices (again from \cite{BGZ}). This result provides the main source of motivation for the definition of companion bases which are the main object of focus of this paper.

In order to state the main result of \cite{BGZ}, we need to introduce some terminology.

\bdfn
\label{quasi cartan}
A symmetrizable matrix $A = (a_{ij})$ is said to be \textit{quasi-Cartan} if $a_{ii} = 2$ for all $i$.
\edfn

\bdfn
\label{positive quasi cartan}
A quasi-Cartan matrix $A$ is said to be \textit{positive} if the symmetrized matrix $DA$ is positive definite.
\edfn

So, we see that a quasi-Cartan matrix $A$ is positive if and only if the principal minors of $A$ are all positive.

The following definition provides a quasi-Cartan analogue of Cartan counterparts.

\bdfn
\label{qcc}
Let $B$ be a skew-symmetrizable matrix. A \textit{quasi-Cartan companion} of $B$ is a quasi-Cartan matrix $A$ such that $|a_{ij}| = |b_{ij}|$ for all $i \neq j$.
\edfn

It is clear that opposite entries have opposite signs in skew-symmetrizable matrices. (Note also that opposite entries have the same sign in symmetrizable matrices.) Therefore, given any skew-symmetrizable matrix $B = (b_{ij})$, we can associate a quiver $\widetilde{\Gamma} (B)$ to $B$ as follows. The vertices of $\widetilde{\Gamma} (B)$ correspond to the rows and columns of $B$, and there is an arrow from the vertex $i$ to the vertex $j$ whenever $b_{ij} >0$.

\bdfn
\label{chordless cycle}
We define a \textit{chordless cycle} in $\widetilde{\Gamma} (B)$ to be a (not necessarily oriented) cycle in $\widetilde{\Gamma} (B)$ such that the full subquiver on its vertices is also a cycle in $\widetilde{\Gamma} (B)$.
\edfn

We can now state the main result from \cite{BGZ} on recognising cluster algebras of finite type (\cite[Theorem 1.2]{BGZ}). We note that an alternative method for recognising cluster algebras of finite type was given in \cite{Sev}.

\bthm
\label{bgz main result}
Let $B$ be a skew-symmetrizable matrix. Then, the cluster algebra associated to (the mutation equivalence class of) $B$ is of finite type if and only if\\
(i) every chordless cycle in $\widetilde{\Gamma} (B)$ is cyclically oriented, and\\
(ii) $B$ has a positive quasi-Cartan companion.
\ethm

It is clear that condition (i) is easy to check for a given skew-symmetrizable matrix $B$, however, condition (ii) could be harder to check as $B$ could have many quasi-Cartan companions. In fact, if $B$ has $N$ non-zero above diagonal entries, then there are $2^N$ different quasi-Cartan companions. The results of the following two propositions (\cite[Proposition 1.4]{BGZ} and \cite[Proposition 1.5]{BGZ} respectively) demonstrate the power of the above theorem. Indeed, it turns out that the positivity of only one (carefully chosen) quasi-Cartan companion of $B$ has to be checked.

\bprop
\label{positive condition}
Let $B$ be skew-symmetrizable, and let $A = (a_{ij})$ be a quasi-Cartan companion of $B$. If $A$ is positive, it must satisfy:
$$
\textrm{For all chordless cycles} \; Z \; \textrm{in} \; \widetilde{\Gamma} (B), \prod_{i \rightarrow j \textrm{ in } Z} (-a_{ij}) < 0. \qquad (\blacktriangle)
$$
\eprop

Now, when choosing a quasi-Cartan companion $A$ for $B$, for each $b_{ij} \neq 0$ with $i \neq j$, we must either choose $a_{ij},a_{ji}>0$ or $a_{ij},a_{ji}<0$. This can be considered as making a sign choice for each arrow in $\widetilde{\Gamma} (B)$. By Proposition~\ref{positive condition}, in order to have any chance of getting a positive quasi-Cartan companion, these signs must be chosen such that each chordless cycle has an odd number of arrows assigned positive sign. The following proposition tells us that if all chordless cycles in $\widetilde{\Gamma} (B)$ are cyclically oriented, then such a choice of quasi-Cartan companion exists. Furthermore, only the positivity of this companion needs to be checked to determine whether or not $\mathcal{A} (B)$ is of finite type. (This is because performing simultaneous sign changes in the rows and columns of a matrix does not affect the positivity of that matrix.)

\bprop
\label{qccexists}
\noindent Let $B$ be a skew-symmetrizable matrix. If every chordless cycle in $\widetilde{\Gamma} (B)$ is cyclically oriented, then $B$ has a quasi-Cartan companion satisfying $(\blacktriangle)$, unique up to simultaneous sign changes in rows and columns.
\eprop

We now recall a standard notion of equivalence for symmetric quasi-Cartan matrices.

\bdfn
\label{equivalent quasi cartan matrices}
Let $A$ and $A'$ be symmetric quasi-Cartan matrices. If there is some integer matrix $E$ with determinant $\pm 1$ such that $A' = E^T AE$, then we say that $A$ and $A'$ are \textit{equivalent}.
\edfn

It is clear that if $A$ is a symmetric positive quasi-Cartan matrix and $A'$ is a symmetric quasi-Cartan matrix equivalent to $A$, then $A'$ is also positive. A result describing the equivalence classes of positive quasi-Cartan matrices in terms of Dynkin types is given in \cite[Proposition 2.9]{BGZ}. The proof of this result gives such a description in the symmetric case (which is the only case we need here) using the above notion of equivalence. We now state this result.

Let $A = (a_{ij})$ be an $n \times n$ symmetric quasi-Cartan matrix. For each $i$, $1 \leq i \leq n$, define an automorphism $s_i$ of the lattice $\mathbb{Z}^n$ by setting $s_i (e_j) = e_j - a_{ij}e_i$ where $\{ e_1 , \ldots , e_n \}$ is the standard basis in $\mathbb{Z}^n$. Let $W(A) \subseteq GL_n (\mathbb{Z})$ be the group generated by $s_1 , \ldots , s_n$.

\bprop
\label{classification of qc matrices}
The following conditions on a symmetric quasi-Cartan matrix $A$ are equivalent:\\
(i) $A$ is positive.\\
(ii) The group $W(A)$ is finite.\\
(iii) There is a root system $\Phi$ and a linearly independent subset $\{ \beta_1 , \ldots , \beta_n \} \subseteq \Phi$ such that $a_{ij} = ( \beta_i , \beta_j )$ for all $1 \leq i,j \leq n$.\\
(iv) $A$ is equivalent to a Cartan matrix $A^0$ of simply-laced Dynkin type.

Under these conditons, if $\Phi_0 \subseteq \Phi$ is the smallest root subsystem of $\Phi$ that contains the set $\{ \beta_1 , \ldots , \beta_n \}$ in (iii), then the Dynkin type of $\Phi_0$ is the same as the Dynkin type of the matrix $A^0$ in (iv), and it characterises $A$ up to equivalence. Furthermore, $W(A)$ is naturally identified with the Weyl group of $\Phi_0$.
\eprop

\note
As noted in \cite{BGZ}, Proposition~\ref{classification of qc matrices} was essentially also proved in \cite[1.2]{Rin1}, using the language of unit quadratic forms.

We now give further consideration to the positive quasi-Cartan companions of those skew-symmetric matrices which give rise to cluster algebras of simply-laced Dynkin type.

Let $B$ be a skew-symmetric matrix and suppose that the cluster algebra $\mathcal{A} (B)$ is of simply-laced Dynkin type. Then, by Theorem~\ref{bgz main result}, $B$ has a (symmetric) positive quasi-Cartan companion $A$. By Proposition~\ref{qccexists}, we see that all positive quasi-Cartan companions for $B$ can be obtained from $A$ by performing simultaneous sign changes in rows and columns, and thus all positive quasi-Cartan companions for $B$ are equivalent. If $B'$ is mutation equivalent to $B$, then $B'$ is skew-symmetric and must also have a positive quasi-Cartan companion. Applying \cite[Proposition 3.2 and Corollary 3.3]{BGZ}, we see that there is a positive quasi-Cartan companion of $B'$ which is equivalent to $A$, and hence all positive quasi-Cartan companions of $B'$ are equivalent to $A$. In particular, by Theorem~\ref{finite type classification}, there must be some skew-symmetric matrix $B_0$, mutation equivalent to $B$, with Cartan counterpart $A_0 = A(B_0)$ a Cartan matrix of simply-laced Dynkin type. It is clear that $A_0$ is a positive quasi-Cartan companion of $B_0$. Therefore, we see that $A_0$ is equivalent to $A$.

Now, let $\Phi$ be a root system of the same Dynkin type as $A_0$ (which is the Dynkin type of $\mathcal{A} (B)$) in some Euclidean space $V$ with positive definite symmetric bilinear form $( \; , \; )$. We have the following simple corollary of Proposition~\ref{classification of qc matrices} which motivates our definition of companion bases in Section~\ref{companion bases}. (We note that its proof is essentially contained in the proof of \cite[Proposition 2.9]{BGZ}.)

\bcor
\label{motivation result}
Let $A = (a_{ij})$ be a positive quasi-Cartan companion of $B$. Then, there is a subset $\{ \beta_1 , \ldots , \beta_n \} \subseteq \Phi$ which is a $\mathbb{Z}$-basis for $\mathbb{Z} \Phi$ such that $a_{ij} = (\beta_i , \beta_j)$ for all $1 \leq i,j \leq n$.
\ecor

\begin{proof}

By Proposition~\ref{classification of qc matrices}, we have that there is a root system $\Phi'$ (in some Euclidean space $V'$ with positive definite symmetric bilinear form $( \; , \; )'$) and a linearly independent subset $\{ \beta_1 , \ldots , \beta_n \} \subseteq \Phi'$ such that $a_{ij} = (\beta_i , \beta_j )'$ for all $1 \leq i,j \leq n$. Furthermore, if $\Phi_0' \subseteq \Phi'$ is the smallest root subsystem of $\Phi'$ that contains $\{ \beta_1 , \ldots , \beta_n \}$, then the Dynkin type of $\Phi_0'$ is the same as the Dynkin type of $A_0$, and hence $\Phi_0'$ and $\Phi$ are isomorphic root systems.

Let $W_{\Phi'}$ be the Weyl group of $\Phi'$, and suppose that $W$ is the subgroup of $W_{\Phi'}$ generated by $s_{\beta_1}, \ldots, s_{\beta_n}$. It is easily seen that $\Phi_0' = W \{ \beta_1 , \ldots , \beta_n \} \subseteq \Phi'$, and a standard argument may then be used to establish that $\{ \beta_1 , \ldots , \beta_n \}$ is a $\mathbb{Z}$-basis for $\mathbb{Z} \Phi_0'$. The result follows since $\Phi_0'$ is isomorphic to $\Phi$.
\end{proof}


\section{Companion Bases}
\label{companion bases}

Motivated by Corollary~\ref{motivation result}, we introduce the definition of a companion basis for the quiver of a cluster-tilted algebra of simply-laced Dynkin type. We then go on to fully describe the relationship between all of the companion bases for such a quiver.

We start by fixing the set-up we will continue to use throughout, except where explicitly stated otherwise. Let $k$ be an algebraically closed field and $Q$ (with $n$ vertices) an alternating quiver of simply-laced Dynkin type, with underlying graph $\Delta$. Let
$$
\mathcal{C} = \frac{\mathcal{D}^b (kQ\mathrm{-mod})}{F}
$$
be the corresponding cluster category. Let $T$ be a basic cluster-tilting object in $\mathcal{C}$ and $\Lambda = \mathrm{End}_{\mathcal{C}} (T)^{\mathrm{op}}$ the corresponding cluster-tilted algebra. Let $\mathcal{A}$ be the cluster algebra of Dynkin type $\Delta$, and suppose that $(\mathbf{x} , B)$ is the seed in $\mathcal{A}$ corresponding to $T$. By Theorem~\ref{quiver of cta is quiver of seed}, $\Gamma = \Gamma (B)$ is the quiver of $\Lambda$. Let $\Gamma_0$ be the set of vertices of $\Gamma$ and $\Gamma_1$ the set of arrows of $\Gamma$. Let $\Phi \subseteq V$ be the root system of Dynkin type $\Delta$ where $V$ is a Euclidean space with positive definite symmetric bilinear form $( \; , \; )$, and let $\Pi = \{ \alpha_1 , \ldots , \alpha_n \}$ be a simple system of $\Phi$. (Note that we may choose $\Phi$ in such a way that the squared length of each root is 2. In particular, this implies that each root is equal to the corresponding coroot, and if $\alpha , \beta \in \Phi$ are non-proportional roots, then $(\alpha , \beta) \in \{ 0 , \pm 1 \}$.)

\bdfn
\label{companion basis}
We call a subset $\{ \gamma_x \colon x \in \Gamma_0 \} \subseteq \Phi$ a \textit{companion basis} for $\Gamma = \Gamma (B)$ if it satisfies the following properties:\\
(i)  $\{ \gamma_x \colon x \in \Gamma_0 \}$ is a $\mathbb{Z}$-basis for $\mathbb{Z} \Phi$.\\
(ii) The matrix $A = (a_{xy})$ given by $a_{xy} = ( \gamma_x , \gamma_y )$ for all $x,y \in \Gamma_0$ is a positive quasi-Cartan companion of $B$.

In this case, we will also often refer to $\{ \gamma_x \colon x \in \Gamma_0 \}$ as a companion basis for $\Gamma$ giving rise to the positive quasi-Cartan companion $A$ of $B$.
\edfn

We note that since the bilinear form $( \; , \; )$ is positive definite and since any $\mathbb{Z}$-basis for $\mathbb{Z} \Phi$ must also be a basis for $V$, then the matrix of inner products of any $\mathbb{Z}$-basis for $\mathbb{Z} \Phi$ must automatically be a positive quasi-Cartan matrix. We can thus modify the definition of a companion basis for $\Gamma$ by replacing condition (ii) above with:\\
(ii)$'$ The matrix $A = (a_{xy})$ given by $a_{xy} = (\gamma_x , \gamma_y)$ for all $x,y \in \Gamma_0$ is a companion of $B$. That is, $|a_{xy}| = |b_{xy}|$ for all $x \neq y$, $x,y \in \Gamma_0$.

The existence of a companion basis for $\Gamma$ is guaranteed by Corollary~\ref{motivation result}. Fix a companion basis $\{ \gamma_x \colon x \in \Gamma_0 \} \subseteq \Phi$ for $\Gamma$, and let $A = (a_{xy})$ be the matrix given by $a_{xy} = ( \gamma_x , \gamma_y)$ for all $x,y \in \Gamma_0$. We can use this companion basis to produce companion bases giving rise to all other positive quasi-Cartan companions of $B$.

Since $\mathcal{A} = \mathcal{A} (B)$ is a cluster algebra of finite type, we see that $A$ is the unique positive quasi-Cartan companion of $B$ up to simultaneous sign changes in rows and columns (due to Propositions \ref{positive condition} and \ref{qccexists}). That is, by applying simultaneous sign changes in the rows and columns of $A$, we can obtain all positive quasi-Cartan companions of $B$. Note also, that all matrices obtained in this way are positive quasi-Cartan companions of $B$.

It is easily checked that if we change the sign of some element in our chosen companion basis for $\Gamma$, then we again get a companion basis for $\Gamma$, this time giving rise to the positive quasi-Cartan matrix obtained from $A$ by simultaneously changing signs in the row and column corresponding to the element whose sign we changed. As a consequence of this, we have the following simple result.

\blemma
\label{sign change result}
Let $\overline{A}$ be a positive quasi-Cartan companion of $B$, obtained from $A$ by applying simultaneous sign changes in some collection $I$ of rows and columns. The subset $\{ \bar{\gamma}_x \colon x \in \Gamma_0 \} \subseteq \Phi$ given by
\begin{eqnarray*}
\bar{\gamma}_x = \left \{ \begin{array}{ll}
\gamma_x & \textrm{ if } x \notin I,\\
-\gamma_x & \textrm{ if } x \in I
\end{array} \right.
\end{eqnarray*}
is a companion basis for $\Gamma$ giving rise to $\overline{A}$.
\elemma

Because an arbitrary companion basis for $\Gamma$ differs from one giving rise to $A$ only by the signs of its elements, we may focus on the relationship between different companion bases for $\Gamma$ giving rise to $A$. Indeed, we will now proceed to show that an arbitrary companion basis for $\Gamma$ giving rise to $A$ can be obtained from the fixed initial companion basis $\{ \gamma_x \colon x \in \Gamma_0 \}$ by applying both an element of the Weyl group $W_{\Phi}$ of $\Phi$ and an orthogonal linear transformation of $V$ that permutes $\Pi$ to each of its elements. We note that a complete description of the orthogonal linear transformations of $V$ that permute $\Pi$ is well known (and easily obtained) in each of the simply-laced Dynkin cases. (See \cite[Chapter VI Section 4]{Bou}, for example.) In each case, these may be interpreted as the graph automorphisms of the corresponding Dynkin diagram.

First, we establish the following result.

\bprop
\label{converse of sigma result}
Let $w \in W_{\Phi}$ and let $\sigma$ be an orthogonal linear transformation of $V$ that permutes $\Pi$. Then, the subset $\{ w\sigma \gamma_x \colon x \in \Gamma_0 \} \subseteq \Phi$ is a companion basis for $\Gamma$ giving rise to $A$.
\eprop

\begin{proof}
We start by checking that $\{ w\sigma \gamma_x \colon x \in \Gamma_0 \} \subseteq \Phi$. Given any root $\alpha$, we have from \cite[Theorem 1.5 \& Corollary 1.5]{Hum} that $\alpha$ can be written in the form $\alpha = s_{\alpha_{i_1}} \cdots s_{\alpha_{i_t}} (\alpha_j)$ for some $\alpha_{i_1} , \ldots , \alpha_{i_t} , \alpha_j \in \Pi$, with $t \in \mathbb{N}$ and $1 \leq i_1 , \ldots , i_t , j \leq n$. Since $\sigma$ is an orthogonal linear transformation, it is easily seen that
$$
\sigma \alpha = s_{\sigma (\alpha_{i_1})} \cdots s_{\sigma (\alpha_{i_t})} (\sigma (\alpha_j)),
$$
and therefore, since $\sigma$ permutes $\Pi$, we deduce that $\sigma \alpha \in \Phi$. Consequently, we see that $\{ w\sigma \gamma_x \colon x \in \Gamma_0 \} \subseteq \Phi$ and that $\sigma$ permutes the set of roots $\Phi$.

Since $w$ and $\sigma$ are both orthogonal linear transformations, it is clear that \linebreak $(w\sigma \gamma_x , w\sigma \gamma_y) = (\gamma_x , \gamma_y) = a_{xy}$ for all $x,y \in \Gamma_0$. The proof is completed by noting that it is easily checked that $\{ w\sigma \gamma_x \colon x \in \Gamma_0 \}$ is a $\mathbb{Z}$-basis for $\mathbb{Z} \Phi$.
\end{proof}

We now devote the remainder of this section to proving the converse of Proposition~\ref{converse of sigma result}. The following well known result will be an important tool in helping us to do this.

\bprop
\label{phi pi equals w pi}
Suppose $\varphi \colon V \rightarrow V$ is an orthogonal linear transformation which permutes the set of roots $\Phi$. Then, there is some $w \in W_{\Phi}$ such that $\varphi \Pi = w \Pi$.
\eprop

\begin{proof}
Since $\varphi$ permutes the set of roots $\Phi$, it is easy to see that $\varphi \Pi$ is a vector space basis for $V$. If $\alpha \in \Phi$, we can write $\alpha = \sum_{i=1}^n c_i \alpha_i$ with either $c_i \geq 0$ for all $1 \leq i \leq n$, or $c_i \leq 0$ for all $1 \leq i \leq n$. But then, since $\varphi \alpha = \sum_{i=1}^n c_i \varphi (\alpha_i)$, and again using the fact that $\varphi$ permutes the set of roots $\Phi$, we see that each $\alpha \in \Phi$ is a linear combination of $\varphi \Pi$ with all coefficients being non-negative, or all coefficients being non-positive. Therefore, $\varphi \Pi$ is a simple system.

Now, from \cite[Theorem 1.4]{Hum} we have that any two simple systems of $\Phi$ are conjugate under $W_{\Phi}$. Therefore, since $\varphi \Pi$ is a simple system, we see that there is some $w \in W_{\Phi}$ such that $\varphi \Pi = w \Pi$.
\end{proof}

Let $\{ \delta_x \colon x \in \Gamma_0 \} \subseteq \Phi$ be another companion basis for $\Gamma$ giving rise to $A$. We aim to describe this companion basis in terms of the companion basis $\{ \gamma_x \colon x \in \Gamma_0 \}$. Define a map $T \colon V \rightarrow V$ by specifying $T (\gamma_x) = \delta_x$ for all $x \in \Gamma_0$, and extending linearly. By definition, we have that $T$ is an invertible linear transformation. Moreover, it is clear that $T$ is an orthogonal transformation. So, in order to be able to apply Proposition~\ref{phi pi equals w pi} to $T$, we must show that $T$ permutes the set of roots $\Phi$. The following result provides the main step towards establishing this. The proof we present is an analogue of the first two parts of the proof of \cite[Theorem 1.5]{Hum}. Before stating this result, we first give a preliminary definition.

\bdfn
\label{beta height}
Let $\alpha \in \Phi$. Since $\{ \gamma_x \colon x \in \Gamma_0 \} \subseteq \Phi$ is a companion basis, we can write $\alpha$ uniquely in the form $\alpha = \sum_{x \in \Gamma_0} c_x \gamma_x$ with $c_x \in \mathbb{Z}$ for all $x \in \Gamma_0$. We then define $\sum_{x \in \Gamma_0} |c_x|$ to be the \textit{height} of $\alpha$ with respect to the companion basis $\{ \gamma_x \colon x \in \Gamma_0 \}$.
\edfn

\bprop
\label{key result for T permutes roots}
Let $W_{\Phi}$ be the Weyl group of $\Phi$ and suppose that $W'$ is the subgroup of $W_{\Phi}$ generated by the reflections $s_{\gamma_x}$ for $x \in \Gamma_0$. Then, for any $\alpha \in \Phi$ there exists $w \in W'$ and $x \in \Gamma_0$ such that $\alpha = w \gamma_x$.
\eprop

\begin{proof}
Let $\alpha \in \Phi$ and consider the non-empty subset $W' \alpha \subseteq \Phi$. Let $\delta$ be an element of $W' \alpha$ of minimal height with respect to the companion basis $\{ \gamma_x \colon x \in \Gamma_0 \}$. We claim that $\delta = \pm \gamma_y$ for some $y \in \Gamma_0$.

Since $\{ \gamma_x \colon x \in \Gamma_0 \} \subseteq \Phi$ is a companion basis, we can write $\delta = \sum_{x \in \Gamma_0} c_x \gamma_x$ with $c_x \in \mathbb{Z}$ for all $x \in \Gamma_0$. We have $0 < (\delta , \delta) = \sum_{x \in \Gamma_0} c_x (\delta , \gamma_x)$, and therefore, there must be some $y$ such that $c_y (\delta , \gamma_y) > 0$. If $\delta = \pm \gamma_y$, then we are done. Suppose that this is not the case and consider the root
\begin{eqnarray*}
s_{\gamma_y} (\delta) & = & \delta - (\delta,\gamma_y) \gamma_y\\
& = & \sum_{x \in \Gamma_0} c_x \gamma_x - (\delta,\gamma_y) \gamma_y\\
& = & \sum_{x \neq y} c_x \gamma_x + (c_y - (\delta,\gamma_y)) \gamma_y \in W' \alpha.
\end{eqnarray*}

Now, we have $c_y (\delta,\gamma_y) > 0$. So, there are two possibilities:

(i) $c_y>0$ and $(\delta,\gamma_y)>0$,

(ii) $c_y<0$ and $(\delta,\gamma_y)<0$.

In case (i), we see that $c_y - (\delta,\gamma_y) = c_y - 1 <c_y$, and moreover, that $|c_y - (\delta,\gamma_y)|<|c_y|$. In case (ii), we see that $|c_y - (\delta,\gamma_y)| = |c_y+1|<|c_y|$. Therefore, in either case we see that the height of $s_{\gamma_y} (\delta)$ with respect to the companion basis $\{ \gamma_x \colon x \in \Gamma_0 \}$ is less than the height of $\delta$ with respect to the companion basis $\{ \gamma_x \colon x \in \Gamma_0 \}$. This is a contradiction. So, we must have that $\delta = \pm \gamma_y$.

In particular, there is some $w \in W'$ such that either $w \alpha = \gamma_y$ or $w \alpha = -\gamma_y$. In the former case, we write $\alpha = w' \gamma_y$ by taking $w' = w^{-1} \in W'$, and in the latter case, we write $\alpha = w' \gamma_y$ by taking $w' = (s_{\gamma_y} w)^{-1} = w^{-1} s_{\gamma_y} \in W'$.
\end{proof}

\note
Let $W'$ be as given in Proposition~\ref{key result for T permutes roots}. The third part of the proof of \cite[Theorem 1.5]{Hum} establishes that $W' = W_{\Phi}$.

\bprop
\label{T permutes roots}
$T$ permutes the set of roots $\Phi$.
\eprop

\begin{proof}
Let $\alpha \in \Phi$. In view of Proposition~\ref{key result for T permutes roots}, we can write $\alpha$ in the form $\alpha = s_{\gamma_{x_1}} \cdots s_{\gamma_{x_t}} ( \gamma_y)$ with $x_1 , \ldots , x_t,y \in \Gamma_0$. Furthermore, since $T$ is an orthogonal linear transformation, we have that $T s_{\gamma} = s_{T (\gamma)}T$ for any root $\gamma$. Therefore, we see that
$$
T \alpha = s_{T (\gamma_{x_1})} \cdots s_{T (\gamma_{x_t})} (T (\gamma_y)) = s_{\delta_{x_1}} \cdots s_{\delta_{x_t}} ( \delta_y).
$$
In particular, $T \alpha$ is a root since $\delta_{x_1} , \ldots , \delta_{x_t} , \delta_y$ are all roots, and the result follows since $T$ is invertible.
\end{proof}

The following lemma provides the final step towards the completion of the proof of the converse of Proposition~\ref{converse of sigma result}.

\blemma
\label{U equals w sigma}
Let $U \colon V \rightarrow V$ be an orthogonal linear transformation that permutes $\Phi$. Then there is some $w \in W_{\Phi}$ and orthogonal linear transformation $\sigma \colon V \rightarrow V$ which permutes $\Pi$ such that $U = w \sigma$.
\elemma

\begin{proof}
By Proposition~\ref{phi pi equals w pi}, we have that there is some $w \in W_{\Phi}$ such that $U \Pi = w \Pi$.

Define $\sigma = w^{-1} U \colon V \rightarrow V$. It is immediate that $\sigma$ is an orthogonal linear transformation. Moreover, since $U\Pi = w\Pi$, we see that $w^{-1}U\Pi = \Pi$, and so $\sigma$ permutes the set of simple roots $\Pi$.
\end{proof}

\note
Lemma~\ref{U equals w sigma} is essentially contained in \cite[p.87]{Sam}.

By combining Proposition~\ref{T permutes roots} and Lemma~\ref{U equals w sigma} we have now established the following result, thus completing the proof of the converse of Proposition~\ref{converse of sigma result}.

\bprop
\label{orthogonal transformation permuting simple roots}
There is some $w \in W_{\Phi}$ and some orthogonal linear transformation $\sigma \colon V \rightarrow V$ which permutes $\Pi$ such that $\delta_x = w\sigma \gamma_x$ for all $x \in \Gamma_0$.
\eprop

We have now established the following result.

\bthm
\label{description of beta sets giving rise to A}
Let $\{ \gamma_x \colon x \in \Gamma_0 \}$ be a companion basis for $\Gamma$ giving rise to $A$. Then, the companion bases for $\Gamma$ that give rise to $A$ are precisely the sets of the form $\{ w\sigma \gamma_x \colon x \in \Gamma_0 \}$, where $w \in W_{\Phi}$ and $\sigma$ is an orthogonal linear transformation of $V$ that permutes $\Pi$.
\ethm

In view of Lemma~\ref{sign change result}, we also have the following corollary, giving a complete desription of all of the companion bases for $\Gamma$ in terms of the (arbitrary) initial companion basis $\{ \gamma_x \colon x \in \Gamma_0 \}$.

\bcor
\label{description of beta sets}
Let $\{ \gamma_x \colon x \in \Gamma_0 \}$ be a companion basis for $\Gamma$. Then, the companion bases for $\Gamma$ are precisely the sets of the form $\{ \varepsilon_x w \sigma \gamma_x \colon x \in \Gamma_0 \}$ where $w \in W_{\Phi}$, $\sigma$ is an orthogonal linear transformation of $V$ that permutes $\Pi$, and $\varepsilon_x \in \{ \pm 1 \}$ for all $x \in \Gamma_0$.
\ecor


\section{A Generalisation of Gabriel's Theorem}
\label{A Generalisation of Gabriel's Theorem}

Due to the classification of the cluster algebras of finite type, we have that there is some seed $(\mathbf{x}_0 , B_0 )$ in $\mathcal{A}$ such that the Cartan counterpart $A (B_0)$ is a Cartan matrix of type $\Delta$. Let $T_0$ be the basic cluster-tilting object in $\mathcal{C}$ corresponding to this seed and $\Lambda_0 = \mathrm{End}_{\mathcal{C}} (T_0)^{\mathrm{op}}$ the associated cluster-tilted algebra. It follows from Theorem~\ref{quiver of cta is quiver of seed} that $\Gamma^0 = \Gamma (B_0)$ is the quiver of $\Lambda_0$.

Since the Cartan counterpart $A(B_0)$ is a Cartan matrix of type $\Delta$, we have that the quiver $\Gamma^0$ must be an orientation of $\Delta$. In particular, the graph underlying $\Gamma^0$ is a tree. It therefore follows from \cite[Theorem 4.2]{BMR3} that $\Lambda_0 \cong k\Gamma^0$, and hence $\Lambda_0$ is (isomorphic to) a path algebra of finite representation type.

Now, applying Gabriel's Theorem to $\Lambda_0$ will help us to deduce a little more information about the companion basis $\Pi$ for $\Gamma^0$. Since $\Pi$ is a simple system of $\Phi$, we can write each positive root $\alpha \in \Phi^+$ uniquely as an integral linear combination of $\alpha_1 , \ldots , \alpha_n$, with all coefficients non-negative. In particular, we can associate a vector to each $\alpha \in \Phi^+$ whose components are the coefficients appearing in the expression for $\alpha$ in terms of $\alpha_1 , \ldots , \alpha_n$. Gabriel's Theorem then tells us that the vectors obtained in this way are the dimension vectors of the finitely generated indecomposable $\Lambda_0$-modules.

Since $\Psi = \{ \gamma_x \colon x \in \Gamma_0 \} \subseteq \Phi$ is a companion basis for $\Gamma$, then $\Psi$ is a $\mathbb{Z}$-basis for $\mathbb{Z} \Phi$. So, we can write each root in $\Phi$ uniquely as an integral linear combination of the elements of the companion basis $\Psi$. This enables us to assign a vector to each root, as follows.

\bdfn
\label{vectors assoc pos roots}
Let $\alpha \in \Phi$ and suppose that $\alpha = \sum_{x \in \Gamma_0} c_x \gamma_x$ with $c_x \in \mathbb{Z}$ for all $x \in \Gamma_0$. We define $d_{\alpha}^{\Psi}$ to be the vector $d_{\alpha}^{\Psi} = (|c_x|)_{x \in \Gamma_0}$.
\edfn

Note that for every $\alpha \in \Phi$, the vector $d_{\alpha}^{\Psi}$ associated to $\alpha$ is the same as the vector $d_{-\alpha}^{\Psi}$ associated to $-\alpha$. For this reason, we will usually restrict our attention to the vectors $d_{\alpha}^{\Psi}$ for $\alpha \in \Phi^+$.

Gabriel's Theorem tells us that the vectors $d_{\alpha}^{\Pi}$ for $\alpha \in \Phi^+$ are the dimension vectors of the finitely generated indecomposable $\Lambda_0$-modules. Motivated by this fact, we introduce the following definition.

\bdfn
\label{strong beta set}
We call $\Psi$ a \textit{strong companion basis} for $\Gamma$ if the vectors $d_{\alpha}^{\Psi}$ for $\alpha \in \Phi^+$ are the dimension vectors of the finitely generated indecomposable $\Lambda$-modules.
\edfn

Having already seen how to find all of the companion bases for $\Gamma$, it is natural to ask whether or not we can decide which of these companion bases are strong. For the remainder of this section, we will focus solely on the Dynkin type $A_n$ case, with $n \in \mathbb{N}$ fixed but arbitrary. We then have that $Q$ is an alternating quiver whose underlying graph is the Dynkin diagram of type $A_n$, that the cluster algebra $\mathcal{A}$ is of Dynkin type $A_n$, and that the root system $\Phi$ is of Dynkin type $A_n$. The main result we will establish is the following.

\bthm
\label{main result}
Let $\Lambda = \mathrm{End}_{\mathcal{C}} (T)^{\mathrm{op}}$ be a cluster-tilted algebra of Dynkin type $A_n$, where $T$ is a basic cluster-tilting object in $\mathcal{C}$. Suppose that $(\mathbf{x} , B)$ is the seed corresponding to $T$ in the cluster algebra $\mathcal{A}$, so that $\Gamma = \Gamma (B)$ is the quiver of $\Lambda$. Then, all companion bases for $\Gamma$ are strong. That is, if $\Psi \subseteq \Phi$ is a companion basis for $\Gamma$, then the vectors $d_{\alpha}^{\Psi}$ for $\alpha \in \Phi^+$ are precisely the dimension vectors of the finitely generated indecomposable $\Lambda$-modules.
\ethm

Giving a proof of this result is our main aim of this section. We start by recalling that the quivers of the cluster-tilted algebras of Dynkin type $A_n$ are precisely the quivers arising from triangulations of a regular $(n+3)$-gon (a result of \cite{CCS1} and \cite{CCS2}). We then examine the structure of these quivers.

Let $\mathbb{P}_{n+3}$ be a regular $(n+3)$-gon. We have that $| \Phi_{\geq -1} | = \frac{1}{2} n(n+1)+n = \frac{1}{2} n(n+3)$ which is equal to the number of diagonals of $\mathbb{P}_{n+3}$. Following \cite[Section 3.5]{FZ3}, let the vertices of $\mathbb{P}_{n+3}$ be $P_1 , P_2 , \ldots , P_{n+3}$, labelled in the anticlockwise direction, and identify the almost positive roots with the diagonals of $\mathbb{P}_{n+3}$ as follows. For $1 \leq i \leq \frac{n+1}{2}$, identify $-\alpha_{2i-1} \in \Phi_{\geq -1}$ with the diagonal joining $P_i$ and $P_{n+3-i}$, and for $1 \leq i \leq \frac{n}{2}$, identify $-\alpha_{2i} \in \Phi_{\geq -1}$ with the diagonal joining $P_{i+1}$ and $P_{n+3-i}$. Finally, for each $1 \leq i \leq j \leq n$, identify the positive root $\alpha_i + \alpha_{i+1} + \ldots + \alpha_j$ with the unique diagonal that crosses precisely the diagonals $-\alpha_i , -\alpha_{i+1} , \ldots , -\alpha_j$. (Note that two diagonals are said to cross if they are distinct and have a common interior point.) We then have the following result from \cite[Proposition 3.14]{FZ3}.

\bprop
\label{clusters correspond to triangulations}
Let $\alpha , \beta \in \Phi_{\geq -1}$. Then,
\begin{eqnarray*}
( \alpha || \beta ) = \left \{ \begin{array}{ll}
1 & \textrm{if the diagonals } \alpha \textrm{ and } \beta \, \textrm{ cross,}\\
0 & \textrm{otherwise}.
\end{array} \right.
\end{eqnarray*}

So, compatible sets are collections of mutually non-crossing diagonals. Therefore, there is a one-to-one correspondence between the set of clusters (or equivalently, seeds) of $\mathcal{A}$ and the set of triangulations of $\mathbb{P}_{n+3}$ by non-crossing diagonals.
\eprop

Let $\mathbb{T}$ be the triangulation of $\mathbb{P}_{n+3}$ corresponding to the seed $(\mathbf{x},B)$. We can associate a connected quiver $Q_{\mathbb{T}}$ to $\mathbb{T}$ as in \cite[Section 2.4]{CCS1}. Take the vertices of $Q_{\mathbb{T}}$ to be the midpoints of the diagonals in $\mathbb{T}$. Let $i$ and $j$ be vertices in $Q_{\mathbb{T}}$ lying on diagonals $d_i$ and $d_j$ respectively. Then, there is an arrow from $i$ to $j$ in $Q_{\mathbb{T}}$ if $d_i$ and $d_j$ bound a common triangle (from the triangulation), and the angle of minimal rotation about the common point of $d_i$ and $d_j$ taking the line through $d_i$ to the line through $d_j$ is in the anticlockwise direction. It follows immediately from \cite[Proposition 12.5]{FZ2} that $Q_{\mathbb{T}}$ is precisely the quiver $\Gamma = \Gamma (B)$ associated to $B$ (taking the vertices of $Q_{\mathbb{T}}$ to be indexed by the corresponding cluster variables in $\mathcal{A}$). In view of Theorem~\ref{quiver of cta is quiver of seed}, we then have that $Q_{\mathbb{T}}$ is the quiver of the cluster-tilted algebra $\Lambda = \mathrm{End}_{\mathcal{C}} (T)^{\mathrm{op}}$. Therefore, the quivers associated to the triangulations of $\mathbb{P}_{n+3}$ are precisely the quivers of the cluster-tilted algebras associated to the basic cluster-tilting objects in $\mathcal{C}$.

We may now consider the structure of the quiver $\Gamma$ by examining the (identical) quiver $Q_{\mathbb{T}}$ more closely. Firstly, we note that all of the triangles in the triangulation $\mathbb{T}$ are of the following three types:

(I) Triangles that consist of one diagonal and two boundary edges of $\mathbb{P}_{n+3}$.

(II) Triangles that consist of two diagonals and one boundary edge of $\mathbb{P}_{n+3}$.

(III) Triangles that consist of three diagonals of $\mathbb{P}_{n+3}$.

\note
Since $n \geq 1$, at least one side of any given triangle in $\mathbb{T}$ must be a diagonal of $\mathbb{P}_{n+3}$.

By the definition of $Q_{\mathbb{T}}$, we have that a triangle in $\mathbb{T}$ of type (I) gives rise to a vertex in $Q_{\mathbb{T}}$, a triangle in $\mathbb{T}$ of type (II) gives rise to an arrow between two vertices in $Q_{\mathbb{T}}$, and a triangle in $\mathbb{T}$ of type (III) gives rise to an oriented $3$-cycle in $Q_{\mathbb{T}}$.

Let $x$ be a vertex in $Q_{\mathbb{T}}$, and let $d_x$ be the corresponding diagonal of $\mathbb{P}_{n+3}$ in $\mathbb{T}$. Then, $d_x$ must bound precisely two triangles in $\mathbb{T}$. If $d_x$ bounds two triangles of type (I), then we must have $n=1$, and $x$ is the only vertex of $Q_{\mathbb{T}}$. If $d_x$ bounds a triangle of type (I) and a triangle of type (II), then $x$ has valency one. If $d_x$ bounds a triangle of type (I) and a triangle of type (III), then $x$ lies on a 3-cycle in $Q_{\mathbb{T}}$ and has valency two. If $d_x$ bounds two triangles of type (II), then $x$ has valency two. If $d_x$ bounds a triangle of type (II) and a triangle of type (III), then $x$ is a vertex at which an arrow meets a 3-cycle, and $x$ has valency three. Finally, if $d_x$ bounds two triangles of type (III), then $x$ is a vertex at which two 3-cycles meet, and $x$ has valency four. This covers all possible cases for the vertex $x$.

\blemma
\label{linear sections and cyclically oriented triangles}
In the underlying (unoriented) graph of $Q_{\mathbb{T}}$, the only cycles are 3-cycles arising from triangles in $\mathbb{T}$ of type (III).
\elemma

\begin{proof}
We proceed by induction on $n$. In the initial case when $n=1$, the result is clear. Now, suppose that $n>1$ and that the result holds for all triangulations of any $(k+3)$-gon with $k<n$. Let $d$ be a diagonal in $\mathbb{T}$. The diagonal $d$ divides the polygon $\mathbb{P}_{n+3}$ into two smaller polygons $\mathbb{P}^{d^+}$ and $\mathbb{P}^{d^-}$, with triangulations $\mathbb{T}^{d^+}$ and $\mathbb{T}^{d^-}$ induced from $\mathbb{T}$. By assumption, in the quivers associated to $\mathbb{T}^{d^+}$ and $\mathbb{T}^{d^-}$, the only cycles are 3-cycles arising from triangles of type (III). Also, it's clear that no diagonal in $\mathbb{T}^{d^+}$ can bound the same triangle in $\mathbb{T}$ as a diagonal in $\mathbb{T}^{d^-}$. This completes the proof as it shows that there can be no cycles in the underlying (unoriented) graph of $Q_{\mathbb{T}}$ containing vertices corresponding to diagonals in $\mathbb{T}^{d^+}$ and vertices corresponding to diagonals of $\mathbb{T}^{d^-}$.
\end{proof}

As a result of Lemma~\ref{linear sections and cyclically oriented triangles}, we see that two 3-cycles in $Q_{\mathbb{T}}$ cannot meet in an arrow.

\note
Both the above description of the possibilities at a vertex in $Q_{\mathbb{T}}$ and the result of Lemma~\ref{linear sections and cyclically oriented triangles} were shown independently in \cite{BV} and \cite{Par1}.

Having now understood the structure of the quiver $\Gamma$, let $I_{\Gamma}$ be the admissible ideal of the path algebra $k\Gamma$ generated by the set of all paths in $\Gamma$ consisting of two consecutive arrows in any given 3-cycle. Due to the result of Lemma~\ref{linear sections and cyclically oriented triangles}, we have the following result (\cite[Theorem 4.1]{CCS2}) as an immediate consequence of \cite[Theorem 4.2]{BMR3}.

\bprop
\label{cta is path algebra with relations}
The cluster-tilted algebra $\Lambda$ is isomorphic to $\frac{k\Gamma}{I_{\Gamma}}$.
\eprop

In view of Proposition~\ref{cta is path algebra with relations}, we have that $\Lambda$ is a gentle algebra (see \cite{ABCP}, for example). Moreover, this means that $\Lambda$ is a string algebra (see \cite{BR}), as the gentle algebras form a subclass of the class of string algebras. This is important as a complete description (up to isomorphism) of the finitely generated indecomposable modules over any string algebra was given in \cite{BR}. Indeed, the isomorphism classes of finitely generated indecomposable modules over $\Lambda$ may be described in terms of the strings in $\Lambda$. A \textit{string} in $\Lambda$ is a reduced walk in the quiver $\Gamma$ which avoids the zero relations and hence does not contain two consecutive arrows of any given 3-cycle in $\Gamma$. The length of a string is the number of arrows it contains, and there is a trivial string (that is, a string of length 0) associated to each vertex in $\Gamma$. Since the algebra $\Lambda$ is determined by its quiver $\Gamma$, we note that we may refer to a string in $\Lambda$ as a string in $\Gamma$ without any confusion.

Let $p$ be a string in $\Gamma$. Following \cite{BR}, we associate to $p$ an indecomposable representation $M(p)$ of $\Gamma$ over $k$, which satisfies the relations in $I_{\Gamma}$, as follows. To each vertex of $\Gamma$ lying on $p$ we associate the vector space $k$, and to each vertex of $\Gamma$ not lying on $p$ we associate the zero space. Also, to each arrow of $\Gamma$ lying on $p$ we associate the identity map on $k$, and to each arrow of $\Gamma$ not lying on $p$ we associate the zero map. We refer to (the $\Lambda$-module corresponding to) $M(p)$ as a \textit{string module} and we have that the set of all string modules contains exactly one representative from each isomorphism class of finitely generated indecomposable $\Lambda$-modules. (We note that there are no band modules in this case as there are no cyclic strings in $\Gamma$.)

Therefore, the dimension vectors of the finitely generated indecomposable $\Lambda$-modules may be described in terms of the strings in $\Gamma$. Indeed, if $p$ is a string in $\Gamma$, then the dimension vector of the string module $M(p)$ associated to $p$ has a one in each position corresponding to a vertex of $\Gamma$ lying on $p$ and a zero in each position corresponding to a vertex of $\Gamma$ not lying on $p$.

Our aim is to show that the vectors $d_{\alpha}^{\Psi}$ for $\alpha \in \Phi^+$ are the dimension vectors of the finitely generated indecomposable $\Lambda$-modules. Our strategy will be to describe the vectors $d_{\alpha}^{\Psi}$ for $\alpha \in \Phi^+$ in terms of the strings in $\Gamma$, and check that this description coincides with that of the dimension vectors of the finitely generated indecomposable $\Lambda$-modules. In order to achieve this, we will construct a bijection between the set of strings in $\Gamma$ and the set of positive roots $\Phi^+$. We start by checking that the number of strings in $\Gamma$ is equal to the number of positive roots in $\Phi$. The following simple result will enable us to do this.

\blemma
\label{unique strings}
There is a unique string between any pair $i,j$ of vertices in $\Gamma$.
\elemma

\begin{proof}
The result follows from the structure of the quiver $\Gamma$, as its only (unoriented) cycles are 3-cycles and a string must not contain two consecutive arrows of any given 3-cycle.
\end{proof}

\blemma
\label{number of strings is number of positive roots}
The number of strings in $\Gamma$ is equal to the number of positive roots in $\Phi$.
\elemma

\begin{proof}
The number of (distinct) non-trivial strings in $\Gamma$ is given by $\begin{pmatrix}n \cr 2\end{pmatrix}$, since there is a unique string in $\Gamma$ associated to each pair of distinct vertices of $\Gamma$. Also, there are $n$ trivial strings in $\Gamma$, one associated to each vertex. Therefore, the total number of strings in $\Gamma$ is
$$
\begin{pmatrix}n \cr 2\end{pmatrix} + n = \frac{n!}{2! (n-2)!} + n = \frac{1}{2} n(n-1) + n = \frac{1}{2} n(n+1),
$$
which is equal to the number of positive roots in $\Phi$.
\end{proof}

In order to construct the bijection between the set of strings in $\Gamma$ and the set of positive roots in $\Phi$, we are going to need some way of associating strings in $\Gamma$ to (positive) roots. The following definition provides this.

\bdfn
\label{support}
Let $\alpha \in \Phi$. We can write $\alpha$ uniquely in the form $\alpha = \sum_{x \in \Gamma_0} c_x \gamma_x$ with $c_x \in \mathbb{Z}$ for all $x \in \Gamma_0$. Define $\mathrm{Supp}_{\Psi}(\alpha) = \{ x \in \Gamma_0 \colon c_x \neq 0 \}$ to be the \textit{support} of $\alpha$. If the elements of $\mathrm{Supp}_{\Psi}(\alpha)$ are precisely the vertices of a string $p$ in $\Gamma$, then we say that $\alpha$ has \textit{support} $p$.
\edfn

Since $\Psi$ is a $\mathbb{Z}$-basis for $\mathbb{Z}\Phi$, it's clear that no root can have two different strings as support. So, given any string in $\Gamma$, we will now show that we can exhibit a positive root which has that string as support. This will then establish the desired bijection. We can deal with the trivial strings immediately. Suppose that $p$ is the trivial string on the vertex $z \in \Gamma_0$. Then, $\gamma_z$ is a root with support $p$. Therefore, either $\gamma_z$ or $-\gamma_z$ is a positive root with support $p$.

We must now consider the non-trivial strings in $\Gamma$. Let $p$ be a non-trivial string in $\Gamma$, and suppose that the vertices of $p$, taken in consecutive order from one end of $p$ to the other, are $x_0 , x_1 , \ldots , x_t$ ($t \geq 1$). Supposing that $0 \leq j<i \leq t$, we see that $x_i$ is joined by an arrow in $\Gamma$ to $x_j$ if and only if $j=i-1$. (This follows immediately from the fact that $p$ is a string in $\Gamma$.) Now, $\Psi = \{ \gamma_x \colon x \in \Gamma_0 \} \subseteq \Phi$ is a companion basis for $\Gamma$. Therefore, supposing that $0 \leq j<i \leq t$, we then have that $\left ( \gamma_{x_i} , \gamma_{x_j} \right ) = \pm 1$ if and only if $j = i-1$, and $\left ( \gamma_{x_i} , \gamma_{x_j} \right ) = 0$ otherwise.

The key result enabling us to find a positive root with support $p$ is the following.

\bprop
\label{type A formula for sequence of simple reflections}
If $1 \leq r \leq t$, then
$$
s_{\gamma_{x_r}} \cdots s_{\gamma_{x_1}} (\gamma_{x_0}) = \gamma_{x_0} + \sum_{k=1}^r (-1)^k \left ( \prod_{l=1}^k \left ( \gamma_{x_l} , \gamma_{x_{l-1}} \right ) \right ) \gamma_{x_k}.
$$
\eprop

\begin{proof}
We proceed by induction on $r$.

In the initial case when $r=1$, we have $s_{\gamma_{x_1}} (\gamma_{x_0}) = \gamma_{x_0} - \left ( \gamma_{x_1} , \gamma_{x_0} \right ) \gamma_{x_1}$, so the result holds in this case.

For our induction hypothesis, we suppose that $r \geq 2$ and the result holds for $r-1$. That is,
$$
s_{\gamma_{x_{r-1}}} \cdots s_{\gamma_{x_1}} (\gamma_{x_0}) = \gamma_{x_0} + \sum_{k=1}^{r-1} (-1)^k \left ( \prod_{l=1}^k \left ( \gamma_{x_l} , \gamma_{x_{l-1}} \right ) \right ) \gamma_{x_k}.
$$

Now, $s_{\gamma_{x_r}} ( s_{\gamma_{x_{r-1}}} \cdots s_{\gamma_{x_1}} (\gamma_{x_0})) = s_{\gamma_{x_{r-1}}} \cdots s_{\gamma_{x_1}} (\gamma_{x_0}) - \left ( \gamma_{x_r} , s_{\gamma_{x_{r-1}}} \cdots s_{\gamma_{x_1}} (\gamma_{x_0}) \right ) \gamma_{x_r}$.

We have
\begin{eqnarray*}
\left ( \gamma_{x_r} , s_{\gamma_{x_{r-1}}} \cdots s_{\gamma_{x_1}} (\gamma_{x_0}) \right ) & = & \left ( \gamma_{x_r} , \gamma_{x_0} + \sum_{k=1}^{r-1} (-1)^k \left ( \prod_{l=1}^k \left ( \gamma_{x_l} , \gamma_{x_{l-1}} \right ) \right ) \gamma_{x_k} \right ) \\
& = & \left ( \gamma_{x_r} , \gamma_{x_0} \right ) + \sum_{k=1}^{r-1} (-1)^k \left ( \prod_{l=1}^k \left ( \gamma_{x_l} , \gamma_{x_{l-1}} \right ) \right ) \left ( \gamma_{x_r} , \gamma_{x_k} \right ).
\end{eqnarray*}

Therefore, since $\left ( \gamma_{x_r} , \gamma_{x_0} \right ) = 0$ and since $\left ( \gamma_{x_r} , \gamma_{x_k} \right ) \neq 0$ only for $k = r-1$, we see that
$$
\left ( \gamma_{x_r} , s_{\gamma_{x_{r-1}}} \cdots s_{\gamma_{x_1}} (\gamma_{x_0}) \right ) = (-1)^{r-1} \prod_{l=1}^r \left ( \gamma_{x_l} , \gamma_{x_{l-1}} \right ).
$$

It therefore follows that
$$
s_{\gamma_{x_r}} s_{\gamma_{x_{r-1}}} \cdots s_{\gamma_{x_1}} (\gamma_{x_0}) = \gamma_{x_0} + \sum_{k=1}^r (-1)^k \left ( \prod_{l=1}^k \left ( \gamma_{x_l} , \gamma_{x_{l-1}} \right ) \right ) \gamma_{x_k}
$$
as required.
\end{proof}

By taking $r=t$ in Proposition~\ref{type A formula for sequence of simple reflections}, we have that
$$
s_{\gamma_{x_t}} \cdots s_{\gamma_{x_1}} (\gamma_{x_0}) = \gamma_{x_0} + \sum_{k=1}^t (-1)^k \left ( \prod_{l=1}^k \left ( \gamma_{x_l} , \gamma_{x_{l-1}} \right ) \right ) \gamma_{x_k}.
$$
We claim that the root $s_{\gamma_{x_t}} \cdots s_{\gamma_{x_1}} (\gamma_{x_0})$ has support $p$. Indeed, we have already noted that $\left ( \gamma_{x_l} , \gamma_{x_{l-1}} \right ) = \pm 1$ for all $1 \leq l \leq t$. So, for each $k$, $1 \leq k \leq t$, we have that $(-1)^k \prod_{l=1}^k \left ( \gamma_{x_l} , \gamma_{x_{l-1}} \right ) = \pm 1$. Therefore, $s_{\gamma_{x_t}} \cdots s_{\gamma_{x_1}} (\gamma_{x_0})$ is of the form
$$
s_{\gamma_{x_t}} \cdots s_{\gamma_{x_1}} (\gamma_{x_0}) = \gamma_{x_0} \pm \gamma_{x_1} \pm \cdots \pm \gamma_{x_t},
$$
showing that $s_{\gamma_{x_t}} \cdots s_{\gamma_{x_1}} (\gamma_{x_0})$ has support $p$. We have thus exhibited a positive root with support $p$, for if $s_{\gamma_{x_t}} \cdots s_{\gamma_{x_1}} (\gamma_{x_0})$ is not a positive root, then $-s_{\gamma_{x_t}} \cdots s_{\gamma_{x_1}} (\gamma_{x_0})$ is a positive root with support $p$. We have now proved the following.

\bprop
\label{su paths correspond to pos roots}
There is a bijective correspondence between the set of strings in $\Gamma$ and the set of positive roots $\Phi^+$, with each string in $\Gamma$ corresponding to the unique positive root which has that string as support.
\eprop

Moreover, we have now completed the proof of Theorem~\ref{main result}.

\begin{proof}[Proof of Theorem~\ref{main result}]
For any given string $p$ in $\Gamma$, there is a unique positive root $\alpha$ which has support $p$. Moreover, we have seen that the vector $d_{\alpha}^{\Psi}$ has a one in each position corresponding to a vertex of $\Gamma$ lying on $p$, and zeros everywhere else. This matches the dimension vector of the string module $M(p)$ associated to $p$.
\end{proof}

We conclude this section with an example.

\bex
\label{main result example}
Let $\Gamma$ be the following quiver.
$$
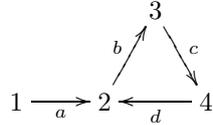
\xymatrix @R=8pt @C=8pt {
 & & & 3 \ar[ddr]^c \\
\\
1 \ar[rr]_a & & 2 \ar[uur]^b & & 4 \ar[ll]^d
}
$$

It is easily checked that $\Gamma$ arises as a quiver associated to a triangulation of a regular 7-gon, and is therefore the quiver of a cluster-tilted algebra $\Lambda$ of Dynkin type $A_4$. In fact, we have $\Lambda \cong \frac{k\Gamma}{I}$ where $I$ is the admissible ideal of $k \Gamma$ given by $I = \langle cb,dc,bd \rangle$ (see Proposition~\ref{cta is path algebra with relations}).

Let $\Phi$ be the root system of Dynkin type $A_4$, and suppose that $\Pi = \{ \alpha_1 , \alpha_2 , \alpha_3 , \alpha_4 \} \subseteq \Phi$ is a simple system of $\Phi$. Then, the set $\Psi = \{ \gamma_1 , \gamma_2 , \gamma_3 , \gamma_4 \} \subseteq \Phi$ given by $\gamma_1 = -\alpha_1$, $\gamma_2 = -\alpha_2 -\alpha_3$, $\gamma_3 = \alpha_3$, $\gamma_4 = \alpha_4$ is a companion basis for $\Gamma$. Indeed, it's clear that $\Psi$ is a $\mathbb{Z}$-basis for $\mathbb{Z} \Phi$. Also, the matrix $A = (a_{ij})$ given by $a_{ij} = ( \gamma_i , \gamma_j )$ for $1 \leq i,j \leq 4$ is a companion (and hence a positive quasi-Cartan companion) of the skew-symmetric integer matrix associated to $\Gamma$.

We find the vectors $d_{\alpha}^{\Psi}$ for $\alpha \in \Phi^+$ by expressing the positive roots in terms of the elements of $\Psi$.
$$
\begin{array}{lll}
\alpha_1 = -\gamma_1 & & d_{\alpha_1}^{\Psi} = (1,0,0,0) \\
\alpha_2 = -\gamma_2 -\gamma_3 & & d_{\alpha_2}^{\Psi} = (0,1,1,0) \\
\alpha_3 = \gamma_3 & & d_{\alpha_3}^{\Psi} = (0,0,1,0) \\
\alpha_4 = \gamma_4 & & d_{\alpha_4}^{\Psi} = (0,0,0,1) \\
\alpha_1 + \alpha_2 = -\gamma_1 -\gamma_2 -\gamma_3 & & d_{\alpha_1 + \alpha_2}^{\Psi} = (1,1,1,0) \\
\alpha_2 + \alpha_3 = -\gamma_2 & & d_{\alpha_2 + \alpha_3}^{\Psi} = (0,1,0,0) \\
\alpha_3 + \alpha_4 = \gamma_3 + \gamma_4 & & d_{\alpha_3 + \alpha_4}^{\Psi} = (0,0,1,1) \\
\alpha_1 + \alpha_2 + \alpha_3 = -\gamma_1 -\gamma_2 & & d_{\alpha_1 +\alpha_2 + \alpha_3}^{\Psi} = (1,1,0,0) \\
\alpha_2 + \alpha_3 + \alpha_4 = -\gamma_2 + \gamma_4 & & d_{\alpha_2 + \alpha_3 + \alpha_4}^{\Psi} = (0,1,0,1) \\
\alpha_1 + \alpha_2 + \alpha_3 + \alpha_4 = -\gamma_1 -\gamma_2 +\gamma_4 & & d_{\alpha_1 + \alpha_2 + \alpha_3 + \alpha_4}^{\Psi} = (1,1,0,1)
\end{array}
$$

These vectors are precisely the dimension vectors of the indecomposable representations of $\Gamma$ that satisfy the relations in $I$, thus confirming that $\Psi$ is strong.
\eex


\section{Companion Basis Mutation and Mutation Maps}
\label{Companion Basis Mutation and Mutation Maps}

We revert to again considering the simply-laced Dynkin case, except where explicitly stated otherwise.

In this section, we establish a companion basis mutation procedure that, given a companion basis for $\Gamma$, produces a companion basis for any mutation of $\Gamma$. Moreover, we note that this provides us with a theoretical, though potentially impractical, method of constructing a companion basis for the quiver of any given cluster-tilted algebra of simply-laced Dynkin type.

We also conjecture that our main result, Theorem~\ref{main result}, may be extended to all of the simply-laced Dynkin cases. This then leads us to consider a collection of maps induced under companion basis mutation.

Let $k$ be a vertex of $\Gamma$. The vertex $k$ corresponds to a row and column of $B$, and hence to a cluster variable $x_k \in \mathbf{x}$. Mutating the seed $(\mathbf{x},B)$ in the direction $x_k$, we obtain a new seed $(\mathbf{x}',B')$. Suppose $T'$ is the basic cluster-tilting object in $\mathcal{C}$ corresponding to the seed $(\mathbf{x}',B')$. Then, we have that $\Gamma' = \Gamma (B')$ is the quiver of the cluster-tilted algebra $\Lambda' = \mathrm{End}_{\mathcal{C}} (T')^{\mathrm{op}}$. We refer to $\Gamma'$ as the quiver obtained from $\Gamma$ by mutating at the vertex $k$.

Although quiver mutation is well understood, it will be useful for our purposes to recall some details. We know (from Remark~\ref{matrix entries remark}) that all matrices appearing in seeds of cluster algebras of Dynkin types $A$, $D$ and $E$ are skew-symmetric integer matrices whose entries belong to the set $\{ 0, \pm 1 \}$. (This means that between any pair of vertices in the quiver of a cluster-tilted algebra of simply-laced Dynkin type there is at most one arrow, and there are no loops on any vertex.) Making use of this information, the following can be obtained via the matrix mutation formula which appears in the definition of a cluster algebra:

Suppose $\Gamma'$ is the quiver obtained from $\Gamma$ by mutating at the vertex $k$. We may identify the vertices of $\Gamma'$ with those of $\Gamma$. Then, the arrows appearing in $\Gamma'$ are the same as those appearing in $\Gamma$, except in the following cases:

(i) All arrows incident with $k$ in $\Gamma$ are replaced in $\Gamma'$ by the corresponding reversed arrows.

(ii) Suppose $x$ and $y$ are vertices such that there is an arrow from $x$ to $k$ in $\Gamma$ and an arrow from $k$ to $y$ in $\Gamma$. Then,

\hspace{15pt} (a) there is an arrow from $x$ to $y$ in $\Gamma'$ if $x$ and $y$ are not joined by an arrow in $\Gamma$,

\hspace{15pt} (b) $x$ and $y$ are not joined by an arrow in $\Gamma'$ if there is an arrow from $y$ to $x$ in $\Gamma$.

The following result introduces the concept of companion basis mutation which gives us a way of obtaining companion bases for $\Gamma'$ (and hence companion bases for any quiver that can be obtained from $\Gamma$ by performing a sequence of quiver mutations) from companion bases for $\Gamma$.

\bthm
\label{beta set mutation theorem}
Let $\{ \gamma_x \colon x \in \Gamma_0 \} \subseteq \Phi$ be a companion basis for $\Gamma$. Then,

(i) the set $\{ \gamma_x' \colon x \in \Gamma_0' \} \subseteq \Phi$ given by
$$
\gamma_x' = \left \{
\begin{array}{cl}
s_{\gamma_k} (\gamma_x) & \textrm{if there is an arrow from } x \textrm{ to } k \textrm{ in } \Gamma,\\
\gamma_x & \textrm{otherwise}
\end{array}
\right.
$$
is a companion basis for $\Gamma'$;

(ii) the set $\{ \gamma_x'' \colon x \in \Gamma_0' \} \subseteq \Phi$ given by
$$
\gamma_x'' = \left \{
\begin{array}{cl}
s_{\gamma_k} (\gamma_x) & \textrm{if there is an arrow from } k \textrm{ to } x \textrm{ in } \Gamma,\\
\gamma_x & \textrm{otherwise}
\end{array}
\right.
$$
is a companion basis for $\Gamma'$.

We refer to $\{ \gamma_x' \colon x \in \Gamma_0' \}$ as the companion basis for $\Gamma'$ obtained from $\{ \gamma_x \colon x \in \Gamma_0 \}$ by mutating inwardly at $k$, and we refer to $\{ \gamma_x'' \colon x \in \Gamma_0' \}$ as the companion basis for $\Gamma'$ obtained from $\{ \gamma_x \colon x \in \Gamma_0 \}$ by mutating outwardly at $k$.
\ethm

\begin{proof}
We only need prove (i); the proof of (ii) is similar.

We have that $\{ \gamma_x \colon x \in \Gamma_0 \} \subseteq \Phi$ is a $\mathbb{Z}$-basis for $\mathbb{Z} \Phi$, and the matrix $A = (a_{xy})$ given by $a_{xy} = (\gamma_x,\gamma_y)$ for all $x,y \in \Gamma_0$ is a positive quasi-Cartan companion of $B$. We must show that $\{ \gamma_x' \colon x \in \Gamma_0' \} \subseteq \Phi$ is a $\mathbb{Z}$-basis for $\mathbb{Z} \Phi$, and the matrix $A' = (a'_{xy})$ given by $a'_{xy} = (\gamma_x',\gamma_y')$ for all $x,y \in \Gamma_0'$ is a companion of $B'$.

We start by checking that $\{ \gamma_x' \colon x \in \Gamma_0' \}$ is a $\mathbb{Z}$-basis for $\mathbb{Z} \Phi$. Let $z \in \mathbb{Z} \Phi$. Since $\{ \gamma_x \colon x \in \Gamma_0 \}$ is a $\mathbb{Z}$-basis for $\mathbb{Z} \Phi$, we can write $z = \sum_{x \in \Gamma_0} a_x \gamma_x$ with $a_x \in \mathbb{Z}$ for all $x \in \Gamma_0$. Let $x$ be a vertex of $\Gamma$. If there is no arrow in $\Gamma$ from $x$ to $k$, then $\gamma_x' = \gamma_x$. On the other hand, if there is an arrow in $\Gamma$ from $x$ to $k$, then $\gamma_x' = s_{\gamma_k} (\gamma_x) = \gamma_x - \left ( \gamma_x , \gamma_k \right ) \gamma_k$, and hence $\gamma_x = \gamma_x' + \left ( \gamma_x , \gamma_k \right ) \gamma_k'$ (since $\gamma_k' = \gamma_k$). Therefore, we have
\begin{eqnarray*}
z & = & \sum_{x \nrightarrow k} a_x \gamma_x + \sum_{x \rightarrow k} a_x \gamma_x\\
& = & \sum_{x \nrightarrow k} a_x \gamma_x' + \sum_{x \rightarrow k} a_x ( \gamma_x' + \left ( \gamma_x , \gamma_k \right ) \gamma_k')\\
& = & \sum_{x \neq k} a_x \gamma_x' + a_k \gamma_k' + \sum_{x \rightarrow k} a_x \left ( \gamma_x , \gamma_k \right ) \gamma_k'\\
& = & \sum_{x \neq k} a_x \gamma_x' + \left ( a_k + \sum_{x \rightarrow k} a_x \left ( \gamma_x , \gamma_k \right ) \right ) \gamma_k'.
\end{eqnarray*}
thus showing that we can write $z$ as an integral linear combination of the roots $\gamma_x'$ for $x \in \Gamma_0'$. To establish the $\mathbb{Z}$-linear independence of the set $\{ \gamma_x' \colon x \in \Gamma_0' \}$, we note that if $\sum_{x \in \Gamma_0'} c_x \gamma_x' = 0$ with $c_x \in \mathbb{Z}$ for all $x \in \Gamma_0'$, then $\sum_{x \neq k} c_x \gamma_x + \left ( c_k - \sum_{x \rightarrow k} c_x (\gamma_x , \gamma_k) \right ) \gamma_k = 0$. Therefore, we must have $c_x = 0$ for all $x \in \Gamma_0'$, due to the $\mathbb{Z}$-linear independence of the set $\{ \gamma_x \colon x \in \Gamma_0 \}$.

Let $x$ and $y$ be vertices in $\Gamma$, with $x \neq k,y$. In order to check that $A'$ is a companion of $B'$, it is sufficient to check that $(\gamma_x',\gamma_y') = \pm 1$ whenever the vertices $x$ and $y$ are joined by an arrow in $\Gamma'$, and $(\gamma_x',\gamma_y') = 0$ otherwise. (We note that all entries of $B'$ belong to the set $\{ 0 , \pm 1 \}$, and moreover, that $x$ and $y$ are joined by an arrow in $\Gamma'$ if and only if $b_{xy}' = \pm 1$.)

Case 1: We start by considering the case $y = k$. There are three possibilities.

Case 1.1: If $x$ and $k$ are not joined by an arrow in $\Gamma$, then they are not joined by an arrow in $\Gamma'$. Also, we have $\gamma_k' = \gamma_k$ and $\gamma_x' = \gamma_x$, and so $(\gamma_k' , \gamma_x') = (\gamma_k , \gamma_x) = 0$.

Case 1.2: If there is an arrow from $k$ to $x$ in $\Gamma$, then there is an arrow from $x$ to $k$ in $\Gamma'$. Also, we have $\gamma_k' = \gamma_k$ and $\gamma_x' = \gamma_x$, and so $(\gamma_k' , \gamma_x') = (\gamma_k , \gamma_x) = \pm 1$.

Case 1.3: If there is an arrow from $x$ to $k$ in $\Gamma$, then there is an arrow from $k$ to $x$ in $\Gamma'$. Also, we have $\gamma_k' = \gamma_k$ and $\gamma_x' = s_{\gamma_k} ( \gamma_x)$, and so, $(\gamma_k' , \gamma_x') = (\gamma_k , s_{\gamma_k} (\gamma_x)) = (s_{\gamma_k} (\gamma_k) , \gamma_x) = -(\gamma_k , \gamma_x) \in \{ \pm 1 \}$.

We now suppose $y \neq k$. If there is no arrow from $x$ to $k$ in $\Gamma$, and no arrow from $y$ to $k$ in $\Gamma$, we see that $\gamma_x' = \gamma_x$ and $\gamma_y' = \gamma_y$, and so $(\gamma_x' , \gamma_y') = (\gamma_x , \gamma_y)$. Note also that the full subquiver of $\Gamma$ on the vertices $x$ and $y$ must be the same as the full subquiver of $\Gamma'$ on the vertices $x$ and $y$. Therefore, from now on, we will assume (without loss of generality) that $x$ is fixed such that there is an arrow from $x$ to $k$ in $\Gamma$.

Case 2: We will first consider the case where $x$ and $y$ are joined by an arrow in $\Gamma$. By assumption, we have $(\gamma_x , \gamma_y) = \pm 1$. There are two possibilities.

Case 2.1: Suppose $k$ and $y$ are not joined by an arrow in $\Gamma$. Then, the arrow joining $x$ and $y$ in $\Gamma$ also appears in $\Gamma'$.
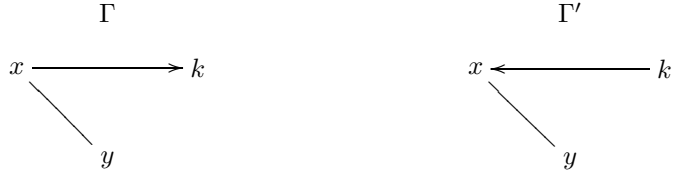
$$
\begin{array}{ccc}
\xymatrix @R=8pt @C=8pt {
 & & {\Gamma} \\
x \ar[rrrr] \ar@{-}[ddrr] & & & & k \\
\\
 & & y
}
&
\hspace{1in}
&
\xymatrix @R=8pt @C=8pt {
 & & {\Gamma'} \\
x \ar@{-}[ddrr] & & & & k \ar[llll] \\
\\
 & & y
}
\end{array}
$$

We have $\gamma_k' = \gamma_k$, $\gamma_x' = s_{\gamma_k} (\gamma_x) = \gamma_x - (\gamma_x , \gamma_k) \gamma_k$, and $\gamma_y' = \gamma_y$. So,
\begin{eqnarray*}
(\gamma_x' , \gamma_y') = (s_{\gamma_k} (\gamma_x) , \gamma_y) & = & (\gamma_x , \gamma_y) - (\gamma_x , \gamma_k)(\gamma_k , \gamma_y)\\
& = & (\gamma_x , \gamma_y)
\end{eqnarray*}
since $(\gamma_k , \gamma_y) = 0$. Hence, $(\gamma_x' , \gamma_y') = \pm 1$.

Case 2.2: Suppose $k$ and $y$ are joined by an arrow in $\Gamma$. Then, there is a triangle in $\Gamma$ on the vertices $k,x,y$. This triangle must be cyclically oriented (since it is a chordless cycle in $\Gamma$). Therefore, since the arrow joining $k$ and $x$ has head $k$, the arrow joining $k$ and $y$ must have head $y$, and the arrow joining $x$ and $y$ must have head $x$. It follows that $x$ and $y$ are not joined by an arrow in $\Gamma'$.
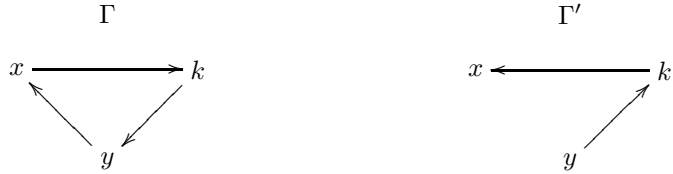
$$
\begin{array}{ccc}
\xymatrix @R=8pt @C=8pt {
 & & {\Gamma} \\
x \ar[rrrr] & & & & k \ar[ddll] \\
\\
 & & y \ar[uull]
}
&
\hspace{1in}
&
\xymatrix @R=8pt @C=8pt {
 & & {\Gamma'} \\
x & & & & k \ar[llll] \\
\\
 & & y \ar[uurr]
}
\end{array}
$$

We have $\gamma_k' = \gamma_k$, $\gamma_x' = s_{\gamma_k} (\gamma_x)$ and $\gamma_y' = \gamma_y$. So, $(\gamma_x' , \gamma_y') = (\gamma_x , \gamma_y) - (\gamma_x , \gamma_k)(\gamma_k , \gamma_y)$.

Now, since $A$ is a positive quasi-Cartan companion of $B$, it follows from Proposition~\ref{positive condition} that an odd number of $(\gamma_x , \gamma_y)$, $(\gamma_x , \gamma_k)$ and $(\gamma_k , \gamma_y)$ are positive. We thus have the following possibilities:
$$
\begin{array}{cccc}
(\gamma_x , \gamma_y) & (\gamma_x , \gamma_k) & (\gamma_k , \gamma_y) & (\gamma_x' , \gamma_y') \\
1 & -1 & -1 & 0 \\
-1 & 1 & -1 & 0 \\
-1 & -1 & 1 & 0 \\
1 & 1 & 1 & 0
\end{array}
$$

Therefore, we see that we must have $(\gamma_x' , \gamma_y') = 0$.

Case 3: We will now consider the case where $x$ and $y$ are not joined by an arrow in $\Gamma$. By assumption, we have $(\gamma_x , \gamma_y) = 0$. There are three possibilities.

Case 3.1: Suppose $k$ and $y$ are not joined by an arrow in $\Gamma$. Then, $x$ and $y$ are not joined by an arrow in $\Gamma'$.
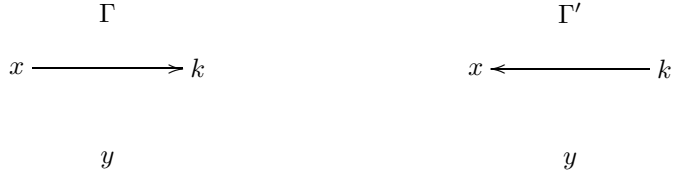
$$
\begin{array}{ccc}
\xymatrix @R=8pt @C=8pt {
 & & {\Gamma} \\
x \ar[rrrr] & & & & k \\
\\
 & & y
}
&
\hspace{1in}
&
\xymatrix @R=8pt @C=8pt {
 & & {\Gamma'} \\
x & & & & k \ar[llll] \\
\\
 & & y
}
\end{array}
$$

We have $\gamma_k' = \gamma_k$, $\gamma_x' = s_{\gamma_k} (\gamma_x)$ and $\gamma_y' = \gamma_y$. So,
\begin{eqnarray*}
(\gamma_x' , \gamma_y') & = & (\gamma_x , \gamma_y) - (\gamma_x , \gamma_k)(\gamma_k , \gamma_y)\\
& = & 0
\end{eqnarray*}
since $(\gamma_x , \gamma_y) = 0$ and $(\gamma_k , \gamma_y) = 0$.

Case 3.2: Suppose that there is an arrow from $y$ to $k$ in $\Gamma$. Again, we see that $x$ and $y$ are not joined by an arrow in $\Gamma'$.
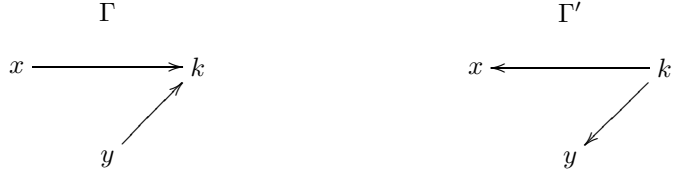
$$
\begin{array}{ccc}
\xymatrix @R=8pt @C=8pt {
 & & {\Gamma} \\
x \ar[rrrr] & & & & k \\
\\
 & & y \ar[uurr]
}
&
\hspace{1in}
&
\xymatrix @R=8pt @C=8pt {
 & & {\Gamma'} \\
x & & & & k \ar[llll] \ar[ddll] \\
\\
 & & y
}
\end{array}
$$

We have $\gamma_k' = \gamma_k$, $\gamma_x' = s_{\gamma_k} (\gamma_x)$ and $\gamma_y' = s_{\gamma_k} (\gamma_y)$. So, $(\gamma_x' , \gamma_y') = (s_{\gamma_k} (\gamma_x) , s_{\gamma_k} (\gamma_y)) = (\gamma_x , \gamma_y ) = 0$.

Case 3.3: Finally, we suppose that there is an arrow from $k$ to $y$ in $\Gamma$. In this case, we see that $k,x,y$ are the vertices of a cyclically oriented triangle in $\Gamma'$. In particular, we have that $x$ and $y$ are joined by an arrow (from $x$ to $y$) in $\Gamma'$.
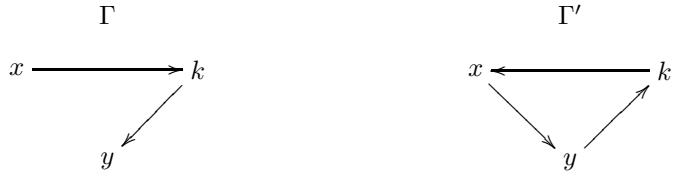
$$
\begin{array}{ccc}
\xymatrix @R=8pt @C=8pt {
 & & {\Gamma} \\
x \ar[rrrr] & & & & k \ar[ddll] \\
\\
 & & y
}
&
\hspace{1in}
&
\xymatrix @R=8pt @C=8pt {
 & & {\Gamma'} \\
x \ar[ddrr] & & & & k \ar[llll] \\
\\
 & & y \ar[uurr]
}
\end{array}
$$

We have $\gamma_k' = \gamma_k$, $\gamma_x' = s_{\gamma_k} (\gamma_x)$ and $\gamma_y' = \gamma_y$. So,
\begin{eqnarray*}
(\gamma_x' , \gamma_y') & = & (\gamma_x , \gamma_y) - (\gamma_x , \gamma_k)(\gamma_k , \gamma_y)\\
& = & - (\gamma_x , \gamma_k)(\gamma_k , \gamma_y) \in \{ \pm 1 \}
\end{eqnarray*}
since $(\gamma_x , \gamma_y) = 0$ and $(\gamma_x , \gamma_k) , (\gamma_k , \gamma_y) \in \{ \pm 1 \}$.

This completes the proof that $\{ \gamma_x' \colon x \in \Gamma_0' \}$ is a companion basis for $\Gamma'$.
\end{proof}

Whilst we have devoted some time to the study of companion bases, we haven't yet discussed how to actually produce a companion basis for $\Gamma$. The above companion basis mutation procedure, at least in principle, provides us with a method for doing this. Recall that due to the classification of the cluster algebras of finite type, we have that there is some seed $(\mathbf{x}_0,B_0)$ in $\mathcal{A}$ such that the Cartan counterpart $A(B_0)$ is a Cartan matrix of type $\Delta$. The seed $(\mathbf{x}_0,B_0)$ corresponds to a basic cluster-tilting object $T_0$ in $\mathcal{C}$, with associated cluster-tilted algebra $\Lambda_0 = \mathrm{End}_{\mathcal{C}} (T_0)^{\mathrm{op}}$. We then have that $\Gamma^0 = \Gamma (B_0)$ is the quiver of $\Lambda_0$, and is an orientation of $\Delta$.

Moreover, recall that the simple system $\Pi = \{ \alpha_1 , \ldots , \alpha_n \}$ is a strong companion basis for $\Gamma^0$ (for some indexing of the simple roots by the vertices of $\Gamma^0$). Since $B$ and $B_0$ are both matrices appearing in seeds in $\mathcal{A}$, then they are mutation equivalent. In particular, we must be able to obtain $\Gamma$ from $\Gamma^0$ by applying a sequence of quiver mutations. So, by applying a corresponding sequence of companion basis mutations to the companion basis $\Pi$ for $\Gamma^0$ (with iterated applications of Theorem~\ref{beta set mutation theorem}), it follows that we obtain a companion basis for $\Gamma$. But of course, this method of finding a companion basis for $\Gamma$ could be difficult to apply, as it requires us to find a sequence of quiver mutations that takes us from $\Gamma^0$ to $\Gamma$.

\note
For more detailed consideration of the construction of explicit companion bases for the quiver of a cluster-tilted algebra, we refer the reader to the forthcoming paper \cite{Par2}. In this paper, we will focus on the cluster-tilted algebras of type $A$ and a special subclass of the cluster-tilted algebras of type $D$, the Nakayama algebras. Indeed, we will give an algorithm that produces a companion basis for the quiver of such a cluster-tilted algebra using only that quiver. Furthermore, we will directly establish that the produced companion bases are strong.

Returning to this paper, we intend to study the procedure of companion basis mutation in more detail. But first, we note the following simple but useful result, showing that the set of vectors associated to the positive roots is the same with respect to any companion basis for $\Gamma$.

\bprop
\label{dalphas do not depend on beta set}
Let $\Psi = \{ \gamma_x \colon x \in \Gamma_0 \} \subseteq \Phi$ and $\Omega = \{ \omega_x \colon x \in \Gamma_0 \} \subseteq \Phi$ be two companion bases for $\Gamma$. Then, $\{ d_{\alpha}^{\Psi} \colon \alpha \in \Phi^+ \} = \{ d_{\alpha}^{\Omega} \colon \alpha \in \Phi^+ \}$. In particular, the set of vectors associated to the set of positive roots, with respect to a given companion basis for $\Gamma$, does not depend on the chosen companion basis.
\eprop

\begin{proof}
Suppose that $\Psi$ gives rise to the positive quasi-Cartan companion $A$ of $B$. By applying sign changes to the elements of $\Omega$, we may obtain a companion basis for $\Gamma$ giving rise to $A$ (see Lemma~\ref{sign change result}). Let $\widetilde{\Omega} = \{ \tilde{\omega}_x \colon x \in \Gamma_0 \} \subseteq \Phi$ be a companion basis for $\Gamma$ giving rise to $A$, where $\tilde{\omega}_x = \varepsilon_x \omega_x$ and $\varepsilon_x \in \{ \pm 1 \}$ for all $x \in \Gamma_0$. It is easily established that $\{ d_{\alpha}^{\Omega} \colon \alpha \in \Phi^+ \} = \{ d_{\alpha}^{\widetilde{\Omega}} \colon \alpha \in \Phi^+ \}$. We may therefore complete the proof by checking that $\{ d_{\alpha}^{\Psi} \colon \alpha \in \Phi^+ \} = \{ d_{\alpha}^{\widetilde{\Omega}} \colon \alpha \in \Phi^+ \}$.

Let $T$ be the (invertible) orthogonal linear transformation of $V$ defined by specifying $T (\gamma_x) = \tilde{\omega}_x$ for all $x \in \Gamma_0$, and extending linearly. We have already seen that $T$ must permute the set of roots $\Phi$ (see Proposition~\ref{T permutes roots}). Let $\alpha \in \Phi$. Then, we can write $\alpha$ uniquely in the form $\alpha = \sum_{x \in \Gamma_0} c_x \gamma_x$ with $c_x \in \mathbb{Z}$ for all $x \in \Gamma_0$. We have $T \alpha = \sum_{x \in \Gamma_0} c_x T (\gamma_x) = \sum_{x \in \Gamma_0} c_x (\tilde{\omega}_x)$, and thus $d_{\alpha}^{\Psi} = d_{T\alpha}^{\widetilde{\Omega}}$. Because $T$ permutes the set of roots $\Phi$, it therefore follows that $\{ d_{\alpha}^{\Psi} \colon \alpha \in \Phi \} = \{ d_{\alpha}^{\widetilde{\Omega}} \colon \alpha \in \Phi \}$, and thus $\{ d_{\alpha}^{\Psi} \colon \alpha \in \Phi^+ \} = \{ d_{\alpha}^{\widetilde{\Omega}} \colon \alpha \in \Phi^+ \}$.
\end{proof}

\bconj
\label{main conjecture}
Let $\Lambda = \mathrm{End}_{\mathcal{C}} (T)^{\mathrm{op}}$ be a cluster-tilted algebra of simply-laced Dynkin type, where $T$ is a basic cluster-tilting object in $\mathcal{C}$. Suppose that $(\mathbf{x} , B)$ is the seed corresponding to $T$ in the cluster algebra $\mathcal{A}$, so that $\Gamma = \Gamma (B)$ is the quiver of $\Lambda$. Then, all companion bases for $\Gamma$ are strong. That is, if $\Psi \subseteq \Phi$ is a companion basis for $\Gamma$, then the vectors $d_{\alpha}^{\Psi}$ for $\alpha \in \Phi^+$ are precisely the dimension vectors of the finitely generated indecomposable $\Lambda$-modules.
\econj

Indeed, we have already established this in the Dynkin type $A$ case. The usefulness of Proposition~\ref{dalphas do not depend on beta set} is that it tells us that in order to prove Conjecture~\ref{main conjecture}, it is enough to establish the existence of a strong companion basis for $\Gamma$. One way of doing this would be to directly construct such a companion basis. (This is the approach taken in \cite{Par2} for the cases mentioned previously.)

A possible alternative approach to establishing the existance of a strong companion basis for $\Gamma$ (thereby proving Conjecture~\ref{main conjecture}) would be to show that companion basis mutation preserves strongness, since we have that $\Pi$ is a strong companion basis for the quiver $\Gamma^0$ of the cluster-tilted algebra $\Lambda_0$.

It therefore makes sense to try to understand the effect of companion basis mutation on the set of vectors associated to the positive roots. To this end, we will now study a collection of maps associated to the companion bases for $\Gamma$, induced by companion basis mutation. In order to enable us to define these maps, we start by showing that the set $\{ d_{\alpha}^{\Psi} \colon \alpha \in \Phi^+ \}$ has $\left | \Phi^+ \right |$ distinct elements. That is, each vector $d_{\alpha}^{\Psi}$ uniquely determines the positive root $\alpha$.

\bprop
\label{vector determines root}
$\left | \{ d_{\alpha}^{\Psi} \colon \alpha \in \Phi^+ \} \right | = \left | \Phi^+ \right |$. That is, the vector $d_{\alpha}^{\Psi}$ is different for each positive root $\alpha$.
\eprop

\begin{proof}
Let $\alpha , \beta \in \Phi^+$ and suppose $d_{\alpha}^{\Psi} = d_{\beta}^{\Psi}$. We will show that $\alpha = \beta$.

Suppose $\alpha = \sum_{x \in \Gamma_0} a_x \gamma_x$ and $\beta = \sum_{x \in \Gamma_0} b_x \gamma_x$ with $a_x , b_x \in \mathbb{Z}$ for all $x \in \Gamma_0$. Since $d_{\alpha}^{\Psi} = d_{\beta}^{\Psi}$, we have $|a_x| = |b_x|$ for all $x \in \Gamma_0$. Let $I = \{ x \colon b_x = a_x \}$ and $J = \{ x \colon b_x \neq a_x \}$ (i.e. $J = \{ x \colon a_x \neq 0 \textrm{ and } b_x = -a_x \}$). Define $\gamma = \sum_{x \in I} a_x \gamma_x \in \mathbb{Z} \Phi$ and $\delta = \sum_{x \in J} a_x \gamma_x \in \mathbb{Z} \Phi$. Then, $\alpha = \gamma + \delta$ and $\beta = \gamma - \delta$.

Since $\alpha , \beta \in \Phi$, we have
$$
2 = ( \alpha , \alpha ) = ( \gamma+\delta , \gamma+\delta ) = (\gamma , \gamma) + 2(\gamma , \delta) + (\delta , \delta)
$$
and
$$
2 = ( \beta , \beta ) = ( \gamma-\delta , \gamma-\delta ) = (\gamma , \gamma) - 2(\gamma , \delta) + (\delta , \delta).
$$
It therefore follows that we must have $(\gamma , \gamma)+(\delta , \delta) = 2$.

Now, for any $0 \neq z \in \mathbb{Z} \Phi$, it is clear that we must have $(z,z) \in \mathbb{N}$. (Note that $( \; , \; )$ is positive definite, and $(z,z)$ is clearly an integer.) But also, we cannot have $(z,z)=1$ (see \cite{CS}, for example). Therefore, there are two cases to consider.

Case 1: Suppose $(\gamma , \gamma) = 2$ and $(\delta , \delta) = 0$. In this case, we must have that $\delta = 0$ and hence $\alpha = \beta$.

Case 2: Suppose $(\gamma , \gamma) = 0$ and $(\delta , \delta) = 2$. In this case, we must have that $\gamma = 0$ and hence $\alpha = -\beta$. But this contradicts the fact that $\alpha$ and $\beta$ are both positive roots.

Therefore, case 2 cannot arise, and so we must have $\alpha = \beta$ as required.
\end{proof}

By Proposition~\ref{dalphas do not depend on beta set}, the set of vectors associated to the positive roots is the same with respect to any companion basis for $\Gamma$. For notational convenience, we will call this set $D (\Gamma)$. In particular, we have $D (\Gamma) = \{ d_{\alpha}^{\Psi} \colon \alpha \in \Phi^+ \}$.

\note
Although we only consider inward companion basis mutation in what follows, we note that analogues of the results established here may be obtained similarly with respect to outward companion basis mutation.

Denote the set of vertices of $\Gamma'$ by $\Gamma_0'$, and let $\Psi' = \{ \gamma_x' \colon x \in \Gamma_0' \} \subseteq \Phi$ be the companion basis for $\Gamma'$ obtained from $\Psi$ by mutating inwardly at $k$. In view of Proposition~\ref{vector determines root}, associated to the companion basis $\Psi$ for $\Gamma$, we have the following bijective map
\begin{eqnarray*}
\phi_{\mathrm{in}}^{\Psi} \colon D (\Gamma) & \longrightarrow & D (\Gamma')\\
d_{\alpha}^{\Psi} & \longmapsto & d_{\alpha}^{\Psi'}.
\end{eqnarray*}
Likewise, we have such a bijective map associated to each companion basis for $\Gamma$. We will show that these bijective maps are all the same. That they all have the same domain, and all have the same codomain, is immediate.

Let $\Omega = \{ \omega_x \colon x \in \Gamma_0 \} \subseteq \Phi$ be another companion basis for $\Gamma$, and suppose that $\Omega' = \{ \omega_x' \colon x \in \Gamma_0' \} \subseteq \Phi$ is the companion basis for $\Gamma'$ obtained from $\Omega$ by mutating inwardly at $k$. We will show that the maps $\phi_{\mathrm{in}}^{\Psi}$ and $\phi_{\mathrm{in}}^{\Omega}$ are equal. By applying sign changes to the elements of $\Omega$, we can obtain a companion basis $\widetilde{\Omega} = \{ \tilde{\omega}_x \colon x \in \Gamma_0 \} \subseteq \Phi$ for $\Gamma$ giving rise to $A$. (Note that $\tilde{\omega}_x = \pm \omega_x$ for all $x \in \Gamma_0$.) Let $\widetilde{\Omega}' = \{ \tilde{\omega}_x' \colon x \in \Gamma_0' \} \subseteq \Phi$ be the companion basis for $\Gamma'$ obtained from $\widetilde{\Omega}$ by mutating inwardly at $k$. It follows as a consequence of the following result (where we consider the situation with a single sign change) that $\phi_{\mathrm{in}}^{\widetilde{\Omega}} = \phi_{\mathrm{in}}^{\Omega}$.

\blemma
\label{mutation map lemma}
Fix $z \in \Gamma_0$ and let $\Theta = \{ \theta_x \colon x \in \Gamma_0 \} \subseteq \Phi$ be the companion basis for $\Gamma$ given by
$$
\theta_x = \left \{ \begin{array}{cl}
-\omega_x & \textrm{if } x = z,\\
\omega_x & \textrm{otherwise}.
\end{array}
\right.
$$
Then: (i) The companion basis $\Theta' = \{ \theta_x' \colon x \in \Gamma_0' \} \subseteq \Phi$ for $\Gamma'$ obtained from $\Theta$ by mutating inwardly at $k$ is given by
$$
\theta_x' = \left \{ \begin{array}{cl}
-\omega_x' & \textrm{if } x = z,\\
\omega_x' & \textrm{otherwise}.
\end{array}
\right.
$$
(ii) $\phi_{\mathrm{in}}^{\Theta} = \phi_{\mathrm{in}}^{\Omega}$.
\elemma

\begin{proof}
(i) There are three cases to consider.

Case 1: Suppose $z \neq k$ and there is no arrow in $\Gamma$ from $z$ to $k$. In this case, we have $\theta_z' = \theta_z = -\omega_z = -\omega_z'$, and $\theta_x' = \omega_x'$ for all $x \neq z$.

Case 2: Suppose there is an arrow in $\Gamma$ from $z$ to $k$. (Note that we must then have $z \neq k$.) In this case, we have $\theta_z' = s_{\theta_k} (\theta_z) = s_{\omega_k} (-\omega_z) = -s_{\omega_k} (\omega_z) = -\omega_z'$, and $\theta_x' = \omega_x'$ for all $x \neq z$.

Case 3: Suppose $z=k$. In this case, we have $\theta_z' = -\omega_z'$ since $\theta_k' = \theta_k = -\omega_k = -\omega_k'$. Also, if $x \neq z$ and there is no arrow in $\Gamma$ from $x$ to $k$, then $\theta_x' = \theta_x = \omega_x = \omega_x'$, and if there is an arrow in $\Gamma$ from $x$ to $k$, then $\theta_x' = s_{\theta_k} (\theta_x) = s_{-\omega_k} (\omega_x) = s_{\omega_k} (\omega_x) = \omega_x'$.

Thus, in each case we see that $\theta_x' = \omega_x'$ for all $x \neq z$, and $\theta_z' = -\omega_z'$.

(ii) Let $\alpha \in \Phi^+$ and write $\alpha = \sum_{x \in \Gamma_0} a_x \omega_x$ with $a_x \in \mathbb{Z}$ for all $x \in \Gamma_0$. Then, we have $\alpha = \sum_{x \in \Gamma_0} b_x \theta_x$ where $b_x = a_x$ for all $x \neq z$, and $b_z = -a_z$. Furthermore, if $\alpha = \sum_{x \in \Gamma_0'} c_x \omega_x'$ with $c_x \in \mathbb{Z}$ for all $x \in \Gamma_0'$, then $\alpha = \sum_{x \in \Gamma_0'} d_x \theta_x'$ where $d_x = c_x$ for all $x \neq z$, and $d_z = -c_z$. In particular, $d_{\alpha}^{\Omega} = d_{\alpha}^{\Theta}$ and $d_{\alpha}^{\Omega'} = d_{\alpha}^{\Theta'}$ for all $\alpha \in \Phi^+$. It follows immediately that the maps $\phi_{\mathrm{in}}^{\Omega}$ and $\phi_{\mathrm{in}}^{\Theta}$ are equal.
\end{proof}

\bcor
\label{mutation map corollary}
$\phi_{\mathrm{in}}^{\Omega} = \phi_{\mathrm{in}}^{\widetilde{\Omega}}$.
\ecor

With the following result, we now prove that all of the bijective maps associated to the companion bases for $\Gamma$ are the same.

\bthm
\label{mutation maps all the same}
The map $\phi_{\mathrm{in}}^{\Psi} \colon D (\Gamma) \rightarrow D (\Gamma')$ does not depend on the companion basis $\Psi$ for $\Gamma$. That is, $\phi_{\mathrm{in}}^{\Psi} = \phi_{\mathrm{in}}^{\Omega}$.
\ethm

\begin{proof}
From Corollary~\ref{mutation map corollary} we have that $\phi_{\mathrm{in}}^{\Omega} = \phi_{\mathrm{in}}^{\widetilde{\Omega}}$. So, we will prove the desired result by showing that $\phi_{\mathrm{in}}^{\Psi} = \phi_{\mathrm{in}}^{\widetilde{\Omega}}$.

Since $\Psi = \{ \gamma_x \colon x \in \Gamma_0 \}$ and $\widetilde{\Omega} = \{ \tilde{\omega}_x \colon x \in \Gamma_0 \}$ are both companion bases for $\Gamma$ giving rise to $A$, we have that the orthogonal linear transformation $T \colon V \rightarrow V$ defined by specifying $T(\gamma_x) = \tilde{\omega}_x$ for all $x \in \Gamma_0$, and extending linearly, permutes the set of roots $\Phi$. Let $\alpha \in \Phi^+$ and write $\alpha = \sum_{x \in \Gamma_0} a_x \gamma_x$ with $a_x \in \mathbb{Z}$ for all $x \in \Gamma_0$. We see that $T\alpha \in \Phi$, and thus either $T\alpha \in \Phi^+$ or $-T\alpha \in \Phi^+$. Suppose without loss of generality that $T\alpha \in \Phi^+$.

We have $T\alpha = \sum_{x \in \Gamma_0} a_x T(\gamma_x) = \sum_{x \in \Gamma_0} a_x \tilde{\omega}_x$, and therefore $d_{\alpha}^{\Psi} = d_{T\alpha}^{\widetilde{\Omega}}$ ( $= d_{-T\alpha}^{\widetilde{\Omega}}$). By Proposition~\ref{vector determines root}, $T\alpha$ must be the unique element of $\Phi^+$ having this property. So, it only remains to check that $d_{\alpha}^{\Psi'} = d_{T\alpha}^{\widetilde{\Omega}'}$. We start by showing that $T(\gamma_x') = \tilde{\omega}_x'$ for all $x \in \Gamma_0'$. There are two cases to consider.

Case 1: Suppose that there is no arrow in $\Gamma$ from $x$ to $k$. Then, $T(\gamma_x') = T(\gamma_x) = \tilde{\omega}_x = \tilde{\omega}_x'$.

Case 2: Suppose that there is an arrow in $\Gamma$ from $x$ to $k$. Then,
\begin{eqnarray*}
T(\gamma_x') & = & T( s_{\gamma_k} (\gamma_x))\\
& = & s_{T(\gamma_k)}(T(\gamma_x))\\
& = & s_{\tilde{\omega}_k} (\tilde{\omega}_x)\\
& = & \tilde{\omega}_x'.
\end{eqnarray*}

Finally, we see that if $\alpha = \sum_{x \in \Gamma_0'} c_x \gamma_x'$ with $c_x \in \mathbb{Z}$ for all $x \in \Gamma_0'$, then $T\alpha = \sum_{x \in \Gamma_0'} c_x T(\gamma_x') = \sum_{x \in \Gamma_0'} c_x \tilde{\omega}_x'$. Therefore, $d_{\alpha}^{\Psi'} = d_{T\alpha}^{\widetilde{\Omega}'}$ ($ = d_{-T\alpha}^{\widetilde{\Omega}'}$) and hence $\phi_{\mathrm{in}}^{\Psi} = \phi_{\mathrm{in}}^{\widetilde{\Omega}}$. This completes the proof.
\end{proof}

Theorem~\ref{mutation maps all the same} tells us that the map $\phi_{\mathrm{in}}^{\Psi} \colon D (\Gamma) \rightarrow D (\Gamma')$ does not depend on the companion basis $\Psi$ for $\Gamma$. That is, if we replace $\Psi$ with any other companion basis for $\Gamma$, then we still get the same map. We will therefore call this map $\phi_{\mathrm{in}}^{\Gamma}$.


\section{Towards a Description of the Mutation Map $\phi_{\mathrm{in}}^{\Gamma}$}
\label{Towards a Description of the Mutation Map}

Naturally, we would like to describe the map $\phi_{\mathrm{in}}^{\Gamma} \colon D (\Gamma) \rightarrow D (\Gamma')$ explicitly. That is, to find a rule that enables us to directly compute the image of any given vector in $D (\Gamma)$ under $\phi_{\mathrm{in}}^{\Gamma}$. We achieve this in the Dynkin type $A_n$ case here, whilst also making some observations that apply to all of the simply-laced Dynkin cases. We start by showing that whenever $\phi_{\mathrm{in}}^{\Gamma}$ (equivalently $\phi_{\mathrm{in}}^{\Psi}$) is applied to a vector in $D (\Gamma)$, then the resultant vector differs from the initial vector in at most one component.

Let $\alpha \in \Phi^+$ and write $\alpha = \sum_{x \in \Gamma_0} a_x \gamma_x$ with $a_x \in \mathbb{Z}$ for all $x \in \Gamma_0$. Using this expression for $\alpha$ in terms of $\Psi$, we can obtain an expression for $\alpha$ in terms of $\Psi'$. Suppose $x$ is a vertex of $\Gamma$. If there is no arrow in $\Gamma$ from $x$ to $k$, then $\gamma_x' = \gamma_x$. On the other hand, if there is an arrow in $\Gamma$ from $x$ to $k$, then $\gamma_x' = s_{\gamma_k} (\gamma_x) = \gamma_x - \left ( \gamma_x , \gamma_k \right ) \gamma_k$, and hence $\gamma_x = \gamma_x' + \left ( \gamma_x , \gamma_k \right ) \gamma_k'$ (since $\gamma_k' = \gamma_k$). It follows, as in the proof of Theorem~\ref{beta set mutation theorem}, that
$$
\alpha = \sum_{x \neq k} a_x \gamma_x' + \left ( a_k + \sum_{x \rightarrow k} a_x \left ( \gamma_x , \gamma_k \right ) \right ) \gamma_k'.
$$

Notice in particular that $d_{\alpha}^{\Psi}$ and $d_{\alpha}^{\Psi'}$ differ in at most one component, the $k$-component. Therefore, in order to describe the map $\phi_{\mathrm{in}}^{\Gamma}$, we should aim to express the $k$-component of $d_{\alpha}^{\Psi'}$ solely in terms of the components of $d_{\alpha}^{\Psi}$. That is, we want to express $\left | a_k + \sum_{x \rightarrow k} a_x \left ( \gamma_x , \gamma_k \right ) \right |$ solely in terms of the integers $|a_x|$ for $x \in \Gamma_0$. The following lemma will subsequently enable us to do this in the Dynkin type $A_n$ case.

\blemma
\label{even difference lemma}
$\left | -|a_k| + \sum_{x \rightarrow k} |a_x| \right | - \left | a_k + \sum_{x \rightarrow k} a_x ( \gamma_x , \gamma_k) \right |$ is an even integer.
\elemma

\begin{proof}
We will show that
$$
-|a_k| + \sum_{x \rightarrow k} |a_x| - \left ( a_k + \sum_{x \rightarrow k} a_x (\gamma_x , \gamma_k) \right ) = -|a_k|-a_k + \sum_{x \rightarrow k} \left ( |a_x| - a_x (\gamma_x , \gamma_k) \right )
$$
is even. The result then follows.

Firstly, since $a_k \in \mathbb{Z}$, we have that $-|a_k|-a_k$ is even.

Let $x$ be a vertex in $\Gamma$, and suppose that there is an arrow in $\Gamma$ from $x$ to $k$. Since $\Psi$ is a companion basis for $\Gamma$, we must have $(\gamma_x , \gamma_k) = \pm 1$. Thus, $|a_x|-a_x ( \gamma_x , \gamma_k) = |a_x| \mp a_x$. Now, since $a_x \in \mathbb{Z}$, $|a_x|-a_x$ and $|a_x|+a_x$ are both even, and hence $|a_x|-a_x ( \gamma_x , \gamma_k)$ is even.

Therefore, $-|a_k|-a_k + \sum_{x \rightarrow k} \left ( |a_x| - a_x (\gamma_x , \gamma_k) \right )$ is even.
\end{proof}

For the remainder of this section we restrict our attention to only the Dynkin type $A_n$ case. (That is, we suppose that $\Delta$ is a Dynkin diagram of type $A_n$.) In this case, we are able to give an explicit description of the map $\phi_{\mathrm{in}}^{\Gamma} \colon D (\Gamma) \rightarrow D (\Gamma')$, by using Lemma~\ref{even difference lemma}. We finish by highlighting a consequence of this description due to Theorem~\ref{main result}.

\bprop
\label{mutation map in type A}
Let $\alpha \in \Phi^+$ and suppose that $\alpha = \sum_{x \in \Gamma_0} a_x \gamma_x$ with $a_x \in \mathbb{Z}$ for all $x \in \Gamma_0$. We then have
$$
\left | a_k + \sum_{x \rightarrow k} a_x (\gamma_x , \gamma_k) \right | = \left | -|a_k| + \sum_{x \rightarrow k} |a_x| \right |.
$$
\eprop

\begin{proof}
We have that all of the components of both $d_{\alpha}^{\Psi}$ and $d_{\alpha}^{\Psi'}$ must belong to the set $\{ 0 , 1 \}$ (see Section~\ref{A Generalisation of Gabriel's Theorem}). Therefore, $|a_x| \in \{ 0,1 \}$ for all vertices $x$ in $\Gamma$, and $\left | a_k + \sum_{x \rightarrow k} a_x (\gamma_x , \gamma_k) \right | \in \{ 0,1 \}$.

Since $|a_x| \in \{ 0,1 \}$ for all vertices $x$ in $\Gamma$, and due to the bijective correspondence between the set of positive roots and the set of strings in $\Gamma$ established in Section~\ref{A Generalisation of Gabriel's Theorem}, a simple case-by-case analysis (based on consideration of the local structure of $\Gamma$ around $k$) establishes that $-|a_k| + \sum_{x \rightarrow k} |a_x| \in \{ 0, \pm 1 \}$.

We have now seen that $\left | a_k + \sum_{x \rightarrow k} a_x (\gamma_x , \gamma_k) \right | , \left | -|a_k| + \sum_{x \rightarrow k} |a_x| \right | \in \{ 0,1 \}$. It therefore follows from Lemma~\ref{even difference lemma} that
$$
\left | a_k + \sum_{x \rightarrow k} a_x (\gamma_x , \gamma_k) \right | = \left | -|a_k| + \sum_{x \rightarrow k} |a_x| \right |
$$
as required.
\end{proof}

We have proved the following.

\bcor
\label{change in dalphas in type A}
Let $\alpha \in \Phi^+$ and suppose that $d_{\alpha}^{\Psi} = (d_x)_{x \in \Gamma_0}$. Then, $d_{\alpha}^{\Psi'}$ is given by
$$
\left ( d_{\alpha}^{\Psi'} \right )_z = \left \{ \begin{array}{ll}
d_z & \textrm{if } z \neq k,\\
\left | -d_k + \sum_{x \rightarrow k} d_x \right | & \textrm{if } z = k.
\end{array} \right.
$$
\ecor

By Theorem~\ref{main result}, we have that the dimension vectors of the finitely generated indecomposable $\Lambda$-modules are precisely the elements of the set $D (\Gamma)$. Likewise, the dimension vectors of the finitely generated indecomposable $\Lambda'$-modules are precisely the elements of the set $D (\Gamma')$. We have therefore established the following.

\bcor
\label{change in dim vectors}
Given the dimension vectors of the finitely generated indecomposable $\Lambda$-modules, we can simply write down the dimension vectors of the finitely generated indecomposable $\Lambda'$-modules.
\ecor


\noindent \textbf{Acknowledgements:}
The research presented in this paper was completed between April 2005 and December 2006 during the author's Ph.D. at the University of Leicester, under the supervision of Robert Marsh, and funded by the Engineering and Physical Sciences Research Council. Consequently, the paper borrows heavily from the author's Ph.D. Thesis \cite{Par1}. The author would like to express his extreme gratitude to Robert Marsh for his unwavering and invaluable support, enthusiasm and insight. For the opportunity to (finally!) write this paper, he would again like to thank Robert Marsh and both the School of Mathematics at the University of Leeds and the Institute for Mathematical Research (FIM, Forschungsinstitut f\"ur Mathematik) at the ETH, Zurich for their kind hospitality during visits in January-February 2011 and April 2011 respectively.

\end{document}